\documentclass[a4paper,11pt]{article}
\usepackage{amsmath}
\usepackage{amsfonts}
\usepackage{amssymb}
\usepackage{graphicx}

\newtheorem{theorem}{Theorem}

\newtheorem{corollary}[theorem]{Corollary}

\newtheorem{definition}[theorem]{Definition}

\newtheorem{lemma}[theorem]{Lemma}

\newtheorem{proposition}[theorem]{Proposition}
\newtheorem{remark}[theorem]{Remark}

\begin{document}

\title{Asymptotics of prediction in functional linear regression with
functional outputs}
\author{ Christophe Crambes and Andr\'{e} Mas\thanks{%
Corresponding author : mas@math.univ-montp2.fr} \\
\\
University of Montpellier}
\date{}
\maketitle

\begin{abstract}
We study prediction in the functional linear model with functional outputs : 
$Y=SX+\varepsilon $ where the covariates $X$ and $Y$ belong to some
functional space and $S$ is a linear operator. We provide the asymptotic
mean square prediction error with exact constants for our estimator which is
based on functional PCA of the input and has a classical form. As a
consequence we derive the optimal choice of the dimension $k_{n}$ of the
projection space. The rates we obtain are optimal in minimax sense and
generalize those found when the output is real. Our main results hold with
no prior assumptions on the rate of decay of the eigenvalues of the input.
This allows to consider a wide class of parameters and inputs $X\left( \cdot
\right) $ that may be either very irregular or very smooth. We also prove a
central limit theorem for the predictor which improves results by Cardot,
Mas and Sarda (2007) in the simpler model with scalar outputs. We show that,
due to the underlying inverse problem, the bare estimate cannot converge in
distribution for the norm of the function space.\bigskip
\end{abstract}


\section{Introduction\protect\bigskip}

\subsection{The model}

Functional data analysis has become these last years an important field in
statistical research, showing a lot of possibilities of applications in many
domains (climatology, teledetection, linguistics, economics, \ldots ). When
one is interested on a phenomenon continuously indexed by time for instance,
it seems appropriate to consider this phenomenon as a whole curve. Practical
aspects also go in this direction, since actual technologies allow to
collect data on thin discretized grids. The papers by Ramsay and Dalzell
(1991) and Frank and Friedman (1993) began to pave the way in favour of this
idea of taking into account the functional nature of these data, and
highlighted the drawbacks of considering a multivariate point of view. Major
references in this domain are the monographs by Ramsay and Silverman (2002,
2005) which give an overview about the philosophy and the basic models
involving functional data. Important nonparametric issues are treated in the
monograph by Ferraty and Vieu (2006).

A particular problem in statistics is to predict the value of an interest
variable $Y$ knowing a covariate $X$. An underlying model can then write :%
\begin{equation*}
Y=r(X)+\varepsilon ,
\end{equation*}

\noindent where $r$ is an operator representing the link between the
variables $X$ and $Y$ and $\varepsilon $ is a noise random variable. In our
functional data context, we want to consider that both variables $X$ and $Y$
are of functional nature, \textit{i.e.} are random functions taking values
on an interval $I=[a,b]$ of $\mathbb{R}$. We assume that $X$ and $Y$ take
values in the space $L^{2}(I)$ of square integrable on $I$. In the following
and in order to simplify, we assume that $I=[0,1]$, which is not restrictive
since the simple transformation $x\longmapsto (x-a)/(b-a)$ allows to come
back to that case.

We assume as well that $X$ and $Y$ are centered. The issue of estimating the
means $\mathbb{E}\left( X\right) $ and $\mathbb{E}\left( Y\right) $ in order
to center the data was exhaustively treated in the literature and is of
minor interest in our setting. The objective of this paper is to consider
the model with functional input and ouptut $:$%
\begin{equation}
Y\left( t\right) =\int_{0}^{1}\mathcal{S}\left( s,t\right) X\left( s\right)
ds+\varepsilon \left( t\right) ,\quad \mathbb{E}\left( \varepsilon |X\right)
=0,  \label{model-kernel}
\end{equation}

\noindent where $\mathcal{S}\left( \cdot ,\cdot \right) $ is an integrable
kernel : $\int \int \left\vert \mathcal{S}\left( s,t\right) \right\vert
dsdt<+\infty $. The kernel $\mathcal{S}$ may be represented on a $3D$-plot
by a surface. The functional historical model (Malfait and Ramsay, 2003) is 
\begin{equation*}
Y\left( t\right) =\int_{0}^{t}\mathcal{S}_{hist}\left( s,t\right) X\left(
s\right) ds+\varepsilon \left( t\right) ,
\end{equation*}%
and may be recovered from the first model be setting $\mathcal{S}\left(
s,t\right) =\mathcal{S}_{hist}\left( s,t\right) 1\!1_{\left\{ s\leq
t\right\} }$ and the surface defining $\mathcal{S}$ is null when $\left(
s,t\right) $ is located in the triangle above the first diagonal of the unit
square.

Model (\ref{model-kernel}) may be viewed as a random Fredholm equation where
both the input an the ouput are random (or noisy). This model has already
been the subject of some studies, as for instance Chiou, M\"{u}ller and Wang
(2004) or Yao, Muller and Wang (2005), which propose an estimation of the
functional parameter $\mathcal{S}$ using functional PCAs of the curves $X$
and $Y$. One of the first studies about this model is due to Cuevas, Febrero
and Fraiman (2002) which considered the case of a fixed design. In this
somewhat different context, they study an estimation of the functional
coefficient of the model and give consistency results for this estimator.
Recently, Antoch \textit{et al.} (2008) proposed a spline estimator of the
functional coefficient in the functional linear model with a functional
response, while Aguilera, Oca\~{n}a and Valderrama (2008) proposed a wavelet
estimation of this coefficient.

We start with a sample $\left( Y_{i}, X_{i} \right)_{1 \leq i \leq n}$ with
the same law as $(Y,X)$, and we consider a new observation $X_{n+1}$. In all
the paper, our goal will be to predict the value of $Y_{n+1}$.

The model (\ref{model-kernel}) may be revisited if one acknowledges that $%
\int_{0}^{1}\mathcal{S}\left( s,t\right) X\left( s\right) ds$ is the image
of $X$ through a general linear integral operator. Denoting $S$ the operator
defined on and with values in $L^{2}\left( \left[ 0,1\right] \right) $ by $%
\left( Sf\right) \left( t\right) =\int_{0}^{1}\mathcal{S}\left( s,t\right)
f\left( s\right) ds$ we obtain from (\ref{model-kernel}) that $Y\left(
t\right) =S\left( X\right) \left( t\right) +\varepsilon \left( t\right) $ or%
\begin{equation*}
Y=SX+\varepsilon ,\quad \text{where} \quad S\left( X\right) \left( t\right) =\int 
\mathcal{S}\left( s,t\right) X\left( s\right) ds.
\end{equation*}%
This fact motivates a more general framework : it may be interesting to
consider Sobolev spaces $W^{m,p}$ instead of $L^{2}\left( \left[ 0,1\right]
\right) $ in order to allow some intrinsic smoothness for the data. It turns
out that, amongst this class of spaces, we should privilege Hilbert spaces.
Indeed the unknown parameter is a linear operator and spectral theory of
these operators acting on Hilbert space allows enough generality, intuitive
approaches and easier practical implementation. That is why in all the
sequel we consider a sample $\left( Y_{i},X_{i}\right) _{1\leq i\leq n}$
where $Y$ and $X$ are independent, identically distributed and take values
in the same Hilbert space $H$ endowed with inner product $\left\langle \cdot
,\cdot \right\rangle $ and associated norm $\left\Vert \cdot \right\Vert .$

Obviously the model we consider generalizes the regression model with a real
output $y$ :%
\begin{equation}
y=\int_{0}^{1}\beta \left( s\right) X\left( s\right) ds+\varepsilon
=\left\langle \beta ,X\right\rangle +\varepsilon ,  \label{scalar-model}
\end{equation}%
and all our results hold in this direction. The literature is wide about (%
\ref{scalar-model}) but we picked articles which are close to our present
concerns and will be cited again later in this work : Yao, M\"{u}ller and
Wang (2005), Hall and Horowitz (2007), Crambes, Kneip, Sarda (2009)...

Since the unknown parameter is here an operator, the infinite-dimensional
equivalent of a matrix, it is worth giving some basic information about
operator theory on Hilbert spaces. The interested reader can find basics and
complements about this topic in the following reference monographs :
Akhiezer and Glazman (1981), Dunford and Schwartz (1988), Gohberg, Goldberg
and Kaashoek (1991). We denote by $\mathcal{L}$ the space of bounded -hence
continuous- operators on a Hilbert space $H$. For our statistical or
probabilistic purposes, we restrain this space to the space of compact
operators $\mathcal{L}_{c}$. Then, any compact and symmetric operator $T$
belonging to $\mathcal{L}_{c}$ admits a unique Schmidt decomposition of the
form $T=\sum_{j\in \mathbb{N}}\mu _{j}\phi _{j}\otimes \phi _{j}$ where the $%
\left( \mu _{j},\phi _{j}\right) $'s are called the eigenelements of $T$,
and the tensor product notation $\otimes $ is defined in the following way:
for any function $f$, $g$ and $h$ belonging to $H$, we define $f\otimes
g=\left\langle g,.\right\rangle f$ or 
\begin{equation*}
\left[ f\otimes g\right] \left( h\right) \left( s\right) =\left( \int
g\left( t\right) h\left( t\right) dt\right) f\left( s\right) .
\end{equation*}

Finally we mention two subclasses of $\mathcal{L}_{c}$ one of which will be
our parameter space. The space of Hilbert-Schmidt operators and trace class
operators are defined respectively by%
\begin{equation*}
\mathcal{L}_{2}=\left\{ T\in \mathcal{L}_{c}:\sum_{j\in \mathbb{N}}\mu
_{j}^{2}<+\infty \right\} ,\ \mathcal{L}_{1}=\left\{ T\in \mathcal{L}%
_{c}:\sum_{j\in \mathbb{N}}\mu _{j}<+\infty \right\} .
\end{equation*}

It is well-known that if $S$ is the linear operator associated to the kernel 
$\mathcal{S}$ like in display (\ref{model-kernel}) then if $\int \int
\left\vert \mathcal{S}\left( s,t\right) \right\vert dsdt<+\infty $, $S$ is
Hilbert-Schmidt and $S$ is trace class if $\mathcal{S}\left( s,t\right) $ is
continuous as a function of $\left( s,t\right) $.

\subsection{Estimation}

Our purpose here is first to introduce the estimator. This estimate looks
basically like the one studied in Yao, M\"{u}ller and Wang (2005). Our
second goal is to justify from a more theoretical position the choice of
such a candidate.

Two strategies may be carried out to propose an estimate of $S.$ They join
finally, like in the finite-dimensional framework. One could consider the
theoretical mean square program (convex in $S$) 
\begin{equation*}
\min_{S\in \mathcal{L}_{2}}\mathbb{E}\left\Vert Y-SX\right\Vert ^{2},
\end{equation*}%
whose solution $S_{\ast }$ is defined by the equation $\mathbb{E}\left[
Y\otimes X\right] =S_{\ast }\mathbb{E}\left[ X\otimes X\right] .$ On the
other hand it is plain that the moment equation :%
\begin{equation*}
\mathbb{E}\left[ Y\otimes X\right] =\mathbb{E}\left[ S\left( X\right)
\otimes X\right] +\mathbb{E}\left[ \varepsilon \otimes X\right]
\end{equation*}%
leads to the same solution. Finally denoting $\Delta =\mathbb{E}\left[
Y\otimes X\right] ,\quad \Gamma =\mathbb{E}\left[ X\otimes X\right] $ we get 
$\Delta =S\Gamma .$ Turning to empirical counterparts with 
\begin{equation*}
\Delta _{n}=\frac{1}{n}\sum_{i=1}^{n}Y_{i}\otimes X_{i},\quad \Gamma _{n}=%
\frac{1}{n}\sum_{i=1}^{n}X_{i}\otimes X_{i},
\end{equation*}%
the estimate $\widehat{S}_{n}$ of $S$ should naturally be defined by $\Delta
_{n}=\widehat{S}_{n}\Gamma _{n}.$Once again the moment method and the
minimization of the mean square program coincide. By the way note that $%
\Delta _{n}=S\Gamma _{n}+U_{n}$ with $U_{n}=\frac{1}{n}\sum_{i=1}^{n}%
\varepsilon _{i}\otimes X_{i}$. The trouble is that, from $\Delta
_{n}=S_{n}\Gamma _{n}$ we cannot directly derive an explicit form for $%
S_{n}. $ Indeed $\Gamma _{n}$ is not invertible on the whole $H$ since it
has finite rank. The next section proposes solutions to solve this inverse
problem by classical methods.


As a last point we note that if $\widehat{S}_{n}$ is an estimate of $S$, a
statistical predictor given a new input $X_{n+1}$ is :%
\begin{equation}
\widehat{Y}_{n+1}\left( t\right) =\widehat{S}_{n}\left( X_{n+1}\right)
\left( t\right) =\int\widehat{\mathcal{S}}\left( s,t\right) X_{n+1}\left(
s\right) ds.  \label{predictor}
\end{equation}

\subsection{Identifiabiliy, inverse problem and regularization issues}

We turn again to the equation which defines the operator $S$ : $\Delta
=S\Gamma .$ Taking a one-to one $\Gamma $ is a first and basic requirement
for identifiability. It is simple to check that if $v\in \ker \Gamma \neq
\left\{ 0\right\} ,$ $\Delta =S\Gamma =\left( S+v\otimes v\right) \Gamma $
for instance and the unicity of $S$ is no more ensured. More precisely, the
inference based on the equation $\Delta =S\Gamma $ does not ensure the
identifiability of the model. From now on we assume that $\ker \Gamma
=\left\{ 0\right\} .$ At this point, some more theoretical concerns should
be mentioned. Indeed, writing $S=\Delta \Gamma ^{-1}$ is untrue. The
operator $\Gamma ^{-1}$ exists whenever $\ker \Gamma =\left\{ 0\right\} $
but is unbounded, that is, not continuous. We refer once again to Dunford
and Schwartz (1988) for instance for developments on unbounded operators. 
It turns out that $\Gamma ^{-1}$ is a linear mapping defined on a dense
domain $\mathcal{D}$ of $H$ which is measurable but continuous at no point
of his domain. Let us denote $\left( \lambda _{j},e_{j}\right) $ the
eigenelements of $\Gamma $. Elementary facts of functional analysis show
that $S_{|\mathcal{D}}=\Delta \Gamma ^{-1}$ where $\mathcal{D}$ is the
domain of $\Gamma ^{-1}$ \textit{i.e.} the range of $\Gamma $ and is defined
by 
\begin{equation*}
\mathcal{D}=\left\{ x=\sum_{j}x_{j}e_{j}\in H:\sum_{j}\frac{x_{j}^{2}}{%
\lambda _{j}^{2}}<+\infty \right\} .
\end{equation*}

A link is possible with probability and gaussian analysis which may be
illustrative. If $\Gamma $ is the covariance operator of a gaussian random
element $X$ on $H$ (a process, a random function, etc) then the Reproducing
Kernel Hilbert Space of $X$ coincides with the domain of $\Gamma ^{-1/2}$
and the range of $\Gamma ^{1/2}$ : $RKHS\left( X\right) =\left\{
x=\sum_{j}x_{j}e_{j}\in H:\sum_{j}x_{j}^{2}/\lambda _{j}<+\infty \right\} .$

The last stumbling stone comes from switching population parameters to
empirical ones. We construct our estimate from the equation $\Delta
_{n}=S\Gamma_{n}+U_{n}$ as seen above and setting $\Delta_{n}=\widehat{S}%
_{n}\Gamma_{n}$. Here the inverse of $\Gamma_{n}$ does not even exist since
this covariance operator is finite-rank. If $\Gamma_{n}$ was invertible we
could set $S_{n}=\Delta_{n}\Gamma_{n}^{-1}$ but we have to regularize $%
\Gamma_{n}$ first. We carry out techniques which are classical in inverse
problems theory. Indeed, the spectral decomposition of $\Gamma_{n}$ is $%
\Gamma_{n}=\sum_{j}\widehat{\lambda}_{j}\left( \widehat{e}_{j}\otimes%
\widehat{e}_{j}\right) $ where $\left( \widehat{\lambda}_{j},\widehat{e}%
_{j}\right) $ are the empirical eigenelements of $\Gamma_{n}$ (the $\widehat{%
\lambda}_{j}$'s are sorted in a decreasing order and some of them may be
null) derived from the functional PCA. The spectral cut regularized inverse
is given for some integer $k$ by 
\begin{equation}
\Gamma_{n}^{\dag}=\sum_{j=1}^{k}\widehat{\lambda}_{j}^{-1}\left( \widehat{e}%
_{j}\otimes\widehat{e}_{j}\right).  \label{gamma-dag}
\end{equation}

The choice of $k=k_{n}$ is crucial ; all the $\left( \widehat{\lambda }%
_{j}\right) _{1\leq j\leq k}$ cannot be null and one should stress that $%
\widehat{\lambda }_{j}^{-1}\uparrow +\infty $ when $j$ increases. The reader
will note that we could define equivalently $\Gamma ^{\dag
}=\sum_{j=1}^{k}\lambda _{j}^{-1}\left( e_{j}\otimes e_{j}\right) .$ From
the definition of the regularized inverse above, we can derive a useful
equation. Indeed, let $\widehat{\Pi }_{k}$ denote the projection of the $k$
first eigenvectors of $\Gamma _{n},$ that is the projection on \textrm{span}$%
\left( \widehat{e}_{1},...,\widehat{e}_{k}\right) .$ Then $\Gamma _{n}^{\dag
}\Gamma _{n}=\Gamma _{n}\Gamma _{n}^{\dag }=\widehat{\Pi }_{k}.$ For further
purpose we define as well $\Pi _{k}$ to be the projection operator on (the
space spanned by) the $k$ first eigenvectors of $\Gamma .$

\begin{remark}
The regularization method we propose is the most intuitive to us but may be
changed by considering : $\Gamma_{n,f}^{\dag}=\sum_{j=1}^{k}f_{n}\left( 
\widehat{\lambda}_{j}\right) \left( \widehat{e}_{j}\otimes\widehat{e}%
_{j}\right) $ where $f_{n}$ is a smooth function which converges pointwise
to $x\rightarrow1/x.$ For instance, we could choose $f_{n}\left( \widehat {%
\lambda}_{j}\right) =\left( \alpha_{n}+\widehat{\lambda}_{j}\right) ^{-1}$
where $\alpha_{n}>0$ and $\alpha_{n}\downarrow0$, and $\Gamma_{n}^{\dag}$
would be the penalized-regularized inverse of $\Gamma_{n}.$ Taking $%
f_{n}\left( \widehat{\lambda}_{j}\right) =\widehat{\lambda}_{j}\left( \alpha
_{n}+\widehat{\lambda}_{j}^{2}\right) ^{-1}$ leads to a Tikhonov
regularization. We refer to the remarks within section 3 of Cardot, Mas,
Sarda (2007) to check that additional assumptions on $f_{n}$ (controlling
the rate of convergence of $f_{n}$ to $x\rightarrow1/x$) allow to generalize
the overall approach of this work to the class of estimates $%
\Gamma_{n,f}^{\dag}$.
\end{remark}

To conclude this subsection, we refer the reader interested by the topic of
inverse problem solving to the following books : Tikhonov and Arsenin
(1977), Groetsch (1993), Engl, Hanke and Neubauer (2000).

\subsection{Assumptions\label{assumptions}}

The assumptions we need are classically of three types : regularity of the
regression parameter $S,$ moment assumptions on $X$ and regularity
assumptions on $X$ which are often expressed in terms of spectral properties
of $\Gamma $ (especially the rate of decrease to zero of its eigenvalues).

\textbf{Assumption on }$S$

As announced sooner, we assume that $S$ is Hilbert Schmidt which may be
rewritten : for any basis $\left( \phi_{j}\right) _{j\in\mathbb{N}}$ of $H$ 
\begin{equation}
\sum_{j,\ell}\left\langle S\left( \phi_{\ell}\right) ,\phi_{j}\right\rangle
^{2}<+\infty.  \label{assumpt-s}
\end{equation}

This assumption finally echoes assumption $\sum_{j}\beta _{j}^{2}<+\infty $
in the functional linear model (\ref{scalar-model}) with real ouptuts. We
already underlined that (\ref{assumpt-s}) is equivalent to assuming that $%
\mathcal{S}$ is doubly integrable if $H$ is $L^{2}\left( \left[ 0,T\right]
\right) $. Finally no continuity or smoothness is required for the kernel $%
\mathcal{S}$ at this point.

\textbf{Moment assumptions on }$X$

In order to better understand the moment assumptions on $X$, we recall the
Karhunen-Loeve development, which is nothing but the decomposition of $X$ in
the basis of the eigenvectors of $\Gamma $, $X=\sum_{j=1}^{+\infty }\sqrt{%
\lambda _{j}}\xi _{j}e_{j}\quad a.s.$ where the $\xi _{j}$'s are independent
centered real random variables with unit variance. We need higher moment
assumptions because we need to apply Bernstein's exponential inequality to
functionals of $\Gamma -\Gamma _{n}.$ We assume that for all $j,\ell \in 
\mathbb{N}$ there exists a constant $b$ such that 
\begin{equation}
\mathbb{E}\left( \left\vert \xi _{j}\right\vert ^{\ell }\right) \leq \frac{%
\ell !}{2}b^{\ell -2}\cdot \mathbb{E}\left( \left\vert \xi _{j}\right\vert
^{2}\right)  \label{assumpt-bernstein}
\end{equation}%
which echoes the assumption (2.19) p. 49 in Bosq (2000). As a consequence,
we see that 
\begin{equation}
\mathbb{E}\left\langle X,e_{j}\right\rangle ^{4}\leq C\left( \mathbb{E}%
\left\langle X,e_{j}\right\rangle ^{2}\right) ^{2}.  \label{H2}
\end{equation}%
This requirement already appears in several papers. It assesses that the
sequence of the fourth moment of the margins of $X$ tends to $0$ quickly
enough. The assumptions above always hold for a gaussian $X$. These
assumptions are close to the moment assumptions usually required when rates
of convergence are addressed.

\textbf{Assumptions on the spectrum of }$\Gamma$

The covariance operator $\Gamma $ is assumed to be injective hence with 
\textit{strictly} positive eigenvalues arranged in a decreasing order. Let
the function $\lambda :\mathbb{R}^{+}\rightarrow \mathbb{R}^{+\ast }$ be
defined by $\lambda \left( j\right) =\lambda _{j}$ for any $j\in \mathbb{N}$
(the $\lambda _{j}$'s are continuously interpolated between $j$ and $j+1.$
>From the assumption above we already know that $\sum_{j}\lambda _{j}<+\infty 
$. Indeed the summability of the eigenvalues of $\Gamma $ is ensured
whenever $\mathbb{E}\left\Vert X\right\Vert ^{2}<+\infty .$ Besides, assume
that for $x$ large enough%
\begin{equation}
x\rightarrow \lambda \left( x\right) \text{ \textrm{is convex}}.  \label{H1}
\end{equation}%
These last conditions are mild and match a very large class of eigenvalues :
with arithmetic decay $\lambda _{j}=Cj^{-1-\alpha }$ where $\alpha >0$ (like
in Hall and Horowitz (2007)), with exponential decay $\lambda
_{j}=Cj^{-\beta }\exp \left( -\alpha j\right) $, Laurent series $\lambda
_{j}=Cj^{-1-\alpha }\left( \log j\right) ^{-\beta }$ or even $\lambda
_{j}=Cj^{-1}\left( \log j\right) ^{-1-\alpha }.$ Such a rate of decay occurs
for extremely irregular processes, even more irregular than the Brownian
motion for which $\lambda _{j}=Cj^{-2}$. In fact our framework initially
relaxes prior assumptions on the rate of decay of the eigenvalues, hence on
the regularity of $X.$ It will be seen later that exact risk and optimality
are obtained when considering specific classes of eigenvalues. Assumption (%
\ref{H1}) is crucial however since the most general Lemmas rely on convex
inequalitites for the eigenvalues.

\section{Asymptotic results}

We are now in a position to introduce our estimate.

\begin{definition}
\label{defestim} The estimate $\widehat{S}_{n}$ of $S$ is defined by : $%
\widehat{S}_{n}=\Delta _{n}\Gamma _{n}^{\dag }$, the associated predictor is 
$\widehat{Y}_{n+1}=\widehat{S}_{n}\left( X_{n+1}\right) =\Delta _{n}\Gamma
_{n}^{\dag }\left( X_{n+1}\right) .$ It is possible to provide a kernel
form. We deduce from $S_{n}=\Delta _{n}\Gamma _{n}^{\dag }$ that 
\begin{equation*}
\mathcal{S}_{n}\left( s,t\right) =\frac{1}{n}\sum_{i=1}^{n}\sum_{j=1}^{k}%
\frac{\int X_{i}\widehat{e}_{j}}{\widehat{\lambda }_{j}}\cdot Y_{i}\left(
t\right) \widehat{e}_{j}\left( s\right) .
\end{equation*}
\end{definition}

Though distinct, this estimate remains close from the one proposed in Yao, M%
\"{u}ller and Wang (2005), the difference consisting in the fact that we do
not consider a Karhunen-Loeve development of $Y$. In the sequel, our main
results are usually given in term of $\widehat{S}_{n}$ but we frequently
switch to the 'kernel' viewpoint since it may be sometimes more
illustrative. Then we implicitely assume that $H=L^{2}\left( \left[ 0,1%
\right] \right) .$

We insist on our philosophy. Estimating $S$ is not our seminal concern. We
focus on the predictor at a random design point $X_{n+1}$, independent from
the initial sample. The issue of estimating $S$ itself may arise typically
for testing. As shown later in this work and as mentioned in Crambes, Kneip
and Sarda (2009), considering the prediction mean square error finally comes
down to studying the mean square error of $S$ for a smooth, intrinsic norm,
depending on $\Gamma$. From now on, all our results are stated when
assumptions of the subsection \ref{assumptions} hold.

\subsection{Mean square prediction error and optimality}

We start with an upper bound from which we deduce, as a Corollary, the exact
asymptotic risk of the predictor. What is considered here is the predictor $%
\widehat{Y}_{n+1}$ based on $\widehat{S}_{n}$ and $X_{n+1}.$ It is compared
with $\mathbb{E}\left( Y_{n+1}|X_{n+1}\right) =S\left( X_{n+1}\right) .$ Let 
$\Gamma_{\varepsilon}=\mathbb{E}\left( \varepsilon\otimes\varepsilon \right) 
$ be the covariance operator of the noise and denote $\sigma
_{\varepsilon}^{2}=\mathrm{tr}\Gamma_{\varepsilon}.$

\begin{theorem}
\label{TH2}The mean square prediction error of our estimate has the
following exact asymptotic development :%
\begin{equation}
\mathbb{E}\left\Vert \widehat{S}_{n}\left( X_{n+1}\right) -S\left(
X_{n+1}\right) \right\Vert ^{2}=\sigma _{\varepsilon }^{2}\frac{k}{n}%
+\sum_{j=k+1}^{+\infty }\lambda _{j}\left\Vert S\left( e_{j}\right)
\right\Vert ^{2}+A_{n}+B_{n},  \label{mse}
\end{equation}%
where $A_{n}\leq C_{A}\left\Vert S\right\Vert _{\mathcal{L}_{2}}k^{2}\lambda
_{k}/n$ and $B_{n}\leq C_{B}k^{2}\log k/n^{2}$ where $C_{A}$ and $C_{B}$ are
constants which do not depend on $k$, $n$ or $S$.
\end{theorem}

The two first term determine the convergence rate : the variance effect
appears through $\sigma _{\varepsilon }^{2}k/n$ and the bias (related to
smoothness) through $\sum_{j=k+1}^{+\infty }\lambda _{j}\left\Vert S\left(
e_{j}\right) \right\Vert ^{2}$. Several comments are needed at this point.
The term $A_{n}$ comes from bias decomposition and $B_{n}$ is a residue from
variance. Both are negligible with respect to the first two terms. Indeed, $%
k\lambda _{k}\rightarrow 0$ since $\sum_{k}\lambda _{k}<+\infty $ and $%
A_{n}=o\left( k/n\right) .$ Turning to $B_{n}$ is a little bit more tricky.
It can be seen from the lines just above the forthcoming Proposition \ref%
{univ-bound} that necessarily $\left( k\log k\right) ^{2}/n\rightarrow 0$
which ensures that $B_{n}=o\left( k/n\right) .$ A second interesting
property arises from Theorem \ref{TH2}. Rewriting $\lambda _{j}\left\Vert
S\left( e_{j}\right) \right\Vert ^{2}=\left\Vert S\Gamma ^{1/2}\left(
e_{j}\right) \right\Vert ^{2}$ we see that the only regularity assumptions
needed may be made from the spectral decomposition of the operator $S\Gamma
^{1/2}$ itself and not from $X$ (or $\Gamma $ as well) and $S$ separately.

Before turning to optimality we introduce the class of parameters $S$ over
which optimality will be obtained.

\begin{definition}
Let $\varphi:\mathbb{R}^{+}\rightarrow\mathbb{R}^{+}$ be a $C^{1}$
decreasing function such that $\sum_{j=1}^{+\infty}\varphi\left( j\right) =1$
and set $\mathcal{L}_{2}\left( \varphi,L\right) $ be the class of linear
operator from $H$ to $H$ be defined by 
\begin{equation*}
\mathcal{L}_{2}\left( \varphi,L\right) =\left\{ T\in\mathcal{L}%
_{2},\left\Vert T\right\Vert _{\mathcal{L}_{2}}\leq L:\left\Vert T\left(
e_{j}\right) \right\Vert \leq L\sqrt{\varphi\left( j\right) }\right\}.
\end{equation*}
\end{definition}

The set $\mathcal{L}_{2}\left( \varphi ,L\right) $ is entirely determined by
the bounding constant $L$ and the function $\varphi $. Horowitz and Hall
(2007) consider the case when $\varphi \left( j\right) =Cj^{-\left( \alpha
+2\beta \right) }$ where $\alpha >1$ and $\beta >1/2.$ As mentioned earlier
we are free here to take any $\varphi $ such that $\int^{+\infty }\varphi
\left( s\right) ds<+\infty $ and which leaves assumption (\ref{H1})
unchanged.

As an easy consequence, we derive the uniform bound with exact constants
below.

\begin{theorem}
Set $L=\left\Vert S\Gamma ^{1/2}\right\Vert _{\mathcal{L}_{2}},$ $\varphi
\left( j\right) =\lambda _{j}\left\Vert S\left( e_{j}\right) \right\Vert
^{2}/L^{2}$ and $k_{n}^{\ast }$ as the integer part of the unique solution
of the integral equation (in $x$) :%
\begin{equation}
\frac{1}{x}\int_{x}^{+\infty }\varphi \left( x\right) dx=\frac{1}{n}\frac{%
\sigma _{\varepsilon }^{2}}{L^{2}}.  \label{k-opt}
\end{equation}%
Let $\mathcal{R}_{n}\left( \varphi ,L\right) $ be the uniform prediction
risk of the estimate $\widehat{S}_{n}$ over the class $\mathcal{L}_{2}\left(
\varphi ,L\right) $ : 
\begin{equation*}
\mathcal{R}_{n}\left( \varphi ,L\right) =\sup_{S\Gamma ^{1/2}\in \mathcal{L}%
_{2}\left( \varphi ,L\right) }\mathbb{E}\left\Vert \widehat{S}_{n}\left(
X_{n+1}\right) -S\left( X_{n+1}\right) \right\Vert ^{2},
\end{equation*}%
then 
\begin{equation*}
\lim \sup_{n\rightarrow +\infty }\frac{n}{k_{n}^{\ast }}\mathcal{R}%
_{n}\left( \varphi ,L\right) =2\sigma _{\varepsilon }^{2}.
\end{equation*}
\end{theorem}

Display (\ref{k-opt}) has a unique solution because the function of $x$ on
the left hand is strictly decreasing. The integer $k_{n}^{\ast }$ is the
optimal dimension : the parameter which minimizes the prediction risk. It
plays the same role as the optimal bandwidth in nonparametric regression.
The upper bound in the display above is obvious from (\ref{mse}). This upper
bound is attained when taking for $S$ the diagonal operator defined in the
basis of eigenvectors by $Se_{j}=L\varphi ^{1/2}\left( j\right) \lambda
_{j}^{-1/2}e_{j}$. The proof of this Theorem is an easy consequence of
Theorem \ref{TH2} hence omitted.

The next Corollary is an attempt to illustrate the consequences of the
previous Theorem by taking explicit sequences $\left( \varphi\left( j\right)
\right) _{j\in\mathbb{N}}$. We chose to treat the case of general Laurent
series (including very irregular input and parameter when $\alpha=0$) and
the case of exponential decay.

\begin{corollary}
Set $\varphi _{a}\left( j\right) =C_{\alpha ,\beta }\left( j^{2+\alpha
}\left( \log j\right) ^{\beta }\right) ^{-1}$ and $\varphi _{b}\left(
j\right) =C_{\alpha }^{\prime }\exp \left( -\alpha j\right) $ where either $%
\alpha >0$ and $\beta \in \mathbb{R}$ or $\alpha =0$ and $\beta >1$, $%
C_{\alpha ,\beta }$ and $C_{\alpha }^{\prime }$ are normalizing constants,
then 
\begin{align*}
\mathcal{R}_{n}\left( \varphi _{a},L\right) & \sim \frac{\left( \log
n\right) ^{\beta /\left( 2+\alpha \right) }}{n^{\left( 1+\alpha \right)
/\left( 2+\alpha \right) }}\left( \frac{C_{\alpha ,\beta }L^{2}}{2\sigma
_{\varepsilon }^{2}}\right) ^{1/\left( 2+\alpha \right) }, \\
\mathcal{R}_{n}\left( \varphi _{b},L\right) & \leq \frac{\log n}{\alpha n}.
\end{align*}
\end{corollary}

In the second display we could not compute an exact bound because equation (%
\ref{k-opt}) has no explicit solution. But the term $\left( \log n\right)
/\alpha n$ is obviously sharp since parametric up to $\log n$. The special
case $\beta =0$ and $\alpha >1$ matches the optimal rate derived in Hall and
Horowitz (2007) with a slight damage due to the fact that the model shows
more complexity ($\mathcal{S}$ is a function of two variables whereas $\beta 
$ the slope parameter in the latter article and in model (\ref{scalar-model}%
) was a function of a single variable). We also refer the reader to Stone
(1982) who underlines this effect of dimension on the convergence rates in
order to check that our result matches the ones announced by Stone.

In our setting the data $Y$ are infinite dimensional. Obtaining lower bound
for optimality in minimax version is slightly different than in the case
studied in Hall and Horowitz (2007), Crambes, Kneip and Sarda (2009). In
order to get a lower bound, our method is close to the one carried out by
Cardot and Johannes (2010), based on a variant of Assouad's Lemma. We
consider gaussian observations under $2^{k_{n}}$ distinct models.

\begin{theorem}
\label{TH2bis}The following bound on the minimax asymptotic risk up to
constants proves that our estimator is optimal in minimax sense :%
\begin{equation*}
\inf_{\widehat{S}_{n}}\sup_{S\in\mathcal{L}_{2}\left( \varphi,L\right) }%
\mathbb{E}\left\Vert \widehat{S}_{n}\left( X_{n+1}\right) -S\left(
X_{n+1}\right) \right\Vert ^{2}\asymp\frac{k_{n}^{\ast}}{n}.
\end{equation*}
\end{theorem}

It appears that another upper bound may be derived from (\ref{mse}). We can
avoid to introduce the class $\mathcal{L}_{2}\left( \varphi ,L\right) .$
>From $\sum_{j}\lambda _{j}=\sigma _{\varepsilon }^{2}$ and $%
\sum_{j}\left\Vert S\left( e_{j}\right) \right\Vert ^{2}=\left\Vert
S\right\Vert _{\mathcal{L}_{2}}^{2}$ we see that the sequences $\lambda _{j}$
and $\left\Vert S\left( e_{j}\right) \right\Vert ^{2}$ may be both bounded
by $j^{-1}(\log j)^{-1}$ hence that $\lambda _{j}\left\Vert S\left(
e_{j}\right) \right\Vert ^{2}\leq j^{-2}(\log j)^{-2}.$ A classical
sum-integral comparison yields then $\sum_{j\geq k+1}\lambda _{j}\left\Vert
S\left( e_{j}\right) \right\Vert ^{2}\leq Ck^{-1}(\log k)^{-2}.$ We obtain
in the Proposition below a new bound for which no regularity assumption is
needed for $S$.

\begin{proposition}
\label{univ-bound}The following bound shows uniformity with respect to all
Hilbert-Schmidt operators $S$ (hence any integrable kernel $\mathcal{S}$)
and all functional data matching the moment assumptions mentioned above :%
\begin{equation*}
\sup_{\left\Vert S\right\Vert _{\mathcal{L}_{2}}\leq L}\mathbb{E}\left\Vert 
\widehat{S}_{n}\left( X_{n+1}\right) -S\left( X_{n+1}\right) \right\Vert
^{2}\leq \sigma _{\varepsilon }^{2}\frac{k}{n}+C\frac{L^{2}}{k\log ^{2}k},
\end{equation*}%
where $C$ is a universal constant. We deduce the uniform bound with no
regularity assumption on the data or on $S$ :%
\begin{equation*}
\lim \sup_{n\rightarrow +\infty }\sqrt{n}\log n\sup_{\left\Vert S\right\Vert
_{\mathcal{L}_{2}}\leq L}\mathbb{E}\left\Vert \widehat{S}_{n}\left(
X_{n+1}\right) -S\left( X_{n+1}\right) \right\Vert ^{2}\leq \sigma
_{\varepsilon }^{2}+CL^{2}.
\end{equation*}
\end{proposition}

The bound above is rough. The constant $C$ does not really matter. The
fundamental idea of the Proposition is to provide an upper bound for the
rate uniformly on balls of $\mathcal{L}_{2}$ without regularity restrictions
: if $\alpha _{n}$ is the rate of prediction error in square norm considered
above, then necessarily $\alpha _{n}\leq n^{-1/2}(\log n)^{-1}$ (in fact we
even have $\alpha _{n}=o\left( n^{-1/2}(\log n)^{-1}\right) $) whatever the
unknown parameter $S$.

\begin{remark}
\label{tradeoff}The bound above holds with highly irregular data (for
instance when $\lambda _{j}\asymp Cj^{-1}(\log j)^{-1-\alpha }$ with $\alpha
>0$ or with very regular data featuring a flat spectrum with $\lambda
_{j}\asymp Cj^{-\gamma }\exp \left( -\alpha j\right) $ or even the
intermediate situation like $\lambda _{j}\asymp Cj^{-1-\beta }(\log
j)^{1+\alpha }$). The literature on linear regression with functional data
usually addressed such issues in restrained case with prior knowledge upon
the eignevalues like $\lambda _{j}\asymp Cj^{-1-\beta }$. The same remarks
are valid when turning to the regularity of the kernel $\mathcal{S}$ or of
the operator $S$ expressed through the sequence $\left\Vert S\left(
e_{j}\right) \right\Vert ^{2}$. Obviously in the case of rapid decay (say at
an exponential rate $\lambda _{j}\asymp C\exp \left( -\alpha j\right) $) one
may argue that multivariate method would fit the data with much accuracy. We
answer that, conversely in such a situation -fitting a linear regression
model- the usual mean square methods turn out to be extremely unstable due
to ill-conditioning. Our method of proof shows that smooth, regular
processes (with rapid decay of $\lambda _{j}$) have good approximation
properties but ill-conditioned $\Gamma _{n}^{\dag }$ (\textit{i.e.} with
rapidly increasing norm) damaging the rate of convergence of $\widehat{S}%
_{n} $ which depends on it. But we readily see that irregular processes
(with slowly decreasing $\lambda _{j}$), despite their poor approximation
properties, lead to a slowly increasing $\Gamma _{n}^{\dag }$ and to solving
an easier inverse problem.
\end{remark}

\begin{remark}
\label{conv-k}At this point it is worth giving a general comment on the rate
of increase of the sequence $k_{n}.$ From the few lines above Proposition %
\ref{univ-bound}, we always have $\left( k\log k\right) ^{2}/n\rightarrow0$
whatever the parameter $S$ in the space of Hilbert-Schmidt operators. This
property will be useful for asymptotics and the mathematical derivations
given in the last section.\bigskip
\end{remark}

\subsection{Weak convergence}

The next and last result deals with weak convergence. We start with a
negative result which shows that due to the underlying inverse problem, the
issue of weak convergence cannot be addressed under too strong topologies.

\begin{theorem}
\label{TH1}It is impossible for $S_{n}$ to converge in distribution for the
Hilbert-Schmidt norm.
\end{theorem}

Once again turning to the predictor, hence smoothing the estimated operator,
will produce a positive result. We improve twofold the results by Cardot,
Mas and Sarda (2007) since first the model is more general and second we
remove the bias term. Weak convergence (convergence in distribution) is
denoted $\overset{w}{\rightarrow }.$ The reader should pay attention to the
fact that the following Theorem holds in space of functions (here $H$).
Within this theorem, two results are proved. The first assesses weak
convergence for the predictor with a bias term. The second removes this bias
at the expense of a more specific assumption on the sequence $k_{n}$.

\begin{theorem}
\label{TH3} If the condition $\left( k\log k\right) ^{2}/n\rightarrow 0$
holds, then 
\begin{equation*}
\sqrt{\frac{n}{k}}\left[ \widehat{S}_{n}\left( X_{n+1}\right) -S\Pi
_{k}\left( X_{n+1}\right) \right] \overset{w}{\rightarrow }\mathcal{G}%
_{\varepsilon }
\end{equation*}%
where $\mathcal{G}_{\varepsilon }$ is a centered gaussian random element
with values in $H$ and covariance operator $\Gamma _{\varepsilon }$.
Besides, denoting $\gamma _{k}=\sup_{j\geq k}\left\{ j\log j\left\Vert
S\left( e_{j}\right) \right\Vert \sqrt{\lambda _{j}}\right\} $ ( it is plain
that $\gamma _{k}\rightarrow 0$) and choosing $k$ such that $n\leq \left(
k\log k\right) ^{2}/\gamma _{k}$ (which means that $\left( k\log k\right)
^{2}/n$ should not decay too quickly to zero), the bias term can be removed
and we obtain 
\begin{equation*}
\sqrt{\frac{n}{k}}\left[ \widehat{S}_{n}\left( X_{n+1}\right) -S\left(
X_{n+1}\right) \right] \overset{w}{\rightarrow }\mathcal{G}_{\varepsilon }.
\end{equation*}
\end{theorem}

\begin{remark}
We pointed out above the improvement in estimating the rate of decrease of
the bias. The proof of the Theorem comes down to proving weak convergence of
a series with values in the space $H$. More precisely, an array $%
\sum_{i=1}^{n}z_{i,n}\varepsilon_{i}$ appears where $z_{i,n}$ are real
valued random variables with increasing variances (when $n\rightarrow+\infty$%
) which are not independent but turn out to be martingale differences.
\end{remark}

>From Theorem \ref{TH3} we deduce general confidence sets for the predictor :
let $\mathcal{K}$ be a continuous set for the measure induced by $\mathcal{G}%
_{\varepsilon }$, that is $\mathbb{P}\left( \mathcal{G}_{\varepsilon }\in
\partial \mathcal{K}\right) =0$ where $\partial \mathcal{K=}\overline{%
\mathcal{K}}\backslash \mathrm{int}\left( \mathcal{K}\right) $ is the
fronteer of $\mathcal{K}$ then $\mathbb{P}\left( \widehat{S}_{n}\left(
X_{n+1}\right) \in S\left( X_{n+1}\right) +\sqrt{\frac{k}{n}}\mathcal{K}%
\right) \rightarrow \mathbb{P}\left( \mathcal{G}_{\varepsilon }\in \mathcal{K%
}\right) $ when $n\rightarrow +\infty $. As an application, we propose the
two following corollaries of Theorem \ref{TH3}. The notation $Y_{n+1}^{\ast
} $ stands for $S\left( X_{n+1}\right) =\mathbb{E}\left(
Y_{n+1}|X_{n+1}\right) $. The first corollary deals with asymptotic
confidence sets for general functionals of the theoretical predictor such as
weighted integrals.

\begin{corollary}
Let $m$ be a fixed function in the space $H=L^{2}\left( \left[ 0,1\right]
\right) $. We have the following asymptotic confidence interval for $\int
Y_{n+1}^{\ast }\left( t\right) m\left( t\right) dt$ at level $1-\alpha $ : 
\begin{equation*}
\mathbb{P}\left( \int_{0}^{1}Y_{n+1}^{\ast }\left( t\right) m\left( t\right)
dt\in \left[ \int_{0}^{1}\widehat{Y}_{n+1}\left( t\right) m\left( t\right)
dt\pm \sqrt{\frac{k}{n}}\sigma _{m}q_{1-\alpha /2}\right] \right) =1-\alpha ,
\end{equation*}%
where $\sigma _{m}^{2}=\left\langle m,\Gamma _{\varepsilon }m\right\rangle
=\int \int \Gamma _{\varepsilon }\left( s,t\right) m\left( t\right) m\left(
s\right) dtds$ rewritten in 'kernel' form and $q_{1-\alpha /2}$ is the
quantile of order $1-\alpha /2$ of the $\mathcal{N}\left( 0,1\right) $
distribution.
\end{corollary}

Theorem \ref{TH3} holds for the Hilbert norm. In order to derive a
confidence interval for $Y_{n+1}^{\ast}\left( t_{0}\right) $ (where $t_{0}$
is fixed in $\left[ 0,1\right] $), we have to make sure that the evaluation
(linear) functional $f\in H\longmapsto f\left( t_{0}\right) $ is continuous
for the norm $\left\Vert \cdot\right\Vert .$ This functional is always
continuous in the space $\left( C\left( \left[ 0,1\right] \right)
,\left\vert \cdot\right\vert _{\infty}\right) $ but is not in the space $%
L^{2}\left( \left[ 0,1\right] \right) .$ A slight change in $H$ will yield
the desired result, stated in the next Corollary.

\begin{corollary}
\label{coro1} When $H=W_{0}^{2,1}\left( \left[ 0,1\right] \right) =\left\{
f\in L^{2}\left( \left[ 0,1\right] \right) :f\left( 0\right) =0,f^{\prime
}\in L^{2}\left( \left[ 0,1\right] \right) \right\} $ endowed with the inner
product $\left\langle u,v\right\rangle =\int_{0}^{1}u^{\prime }v^{\prime }$,
the evaluation functional is continuous with respect to the norm of $H$ and
we can derive from Theorem \ref{TH3} :%
\begin{equation*}
\mathbb{P}\left( Y_{n+1}^{\ast }\left( t_{0}\right) \in \left[ \widehat{Y}%
_{n+1}\left( t_{0}\right) \pm \sqrt{\frac{k}{n}}\sigma _{t_{0}}q_{1-\alpha
/2}\right] \right) =1-\alpha
\end{equation*}%
where $\sigma _{t_{0}}^{2}=\Gamma _{\varepsilon }\left( t_{0},t_{0}\right) .$
\end{corollary}

Note that data $\left( Y_{i}\right) _{1\leq i\leq n}$ reconstructed by cubic
splines and correctly rescaled to match the condition $\left[ f\left(
0\right) =0\right] $ belong to the space $W_{0}^{2,1}\left( \left[ 0,1\right]
\right) $ mentioned in the Corollary.

\begin{remark}
It is out of the scope of this article to go through all the testing issues
which can be solved by Theorem \ref{TH3}. It is interesting to note that if $%
S=0$, the Theorem ensures that%
\begin{equation*}
\sqrt{\frac{n}{k}}\left[ \widehat{S}_{n}\left( X_{n+1}\right) \right] 
\overset{w}{\rightarrow }\mathcal{G}_{\varepsilon },
\end{equation*}%
which may be the starting point for a testing procedure of $S=0$ versus
various alternatives.
\end{remark}

\subsection{Comparison with existing results - Conclusion}

The literature on linear models for functional data gave birth to impressive
and brilliant recent works. We discuss briefly here our contribution with
respect to some articles, close in spirit to this present paper.

We consider exactly the same model (with functional outputs) as Yao, M\"{u}%
ller and Wang (2005) and our estimate is particularly close to the one they
propose. In their work the case of longitudinal data was studied with care
with possibly sparse and irregular data. They introduce a very interesting
functional version of the $R^{2}$ and prove convergence in probability of
their estimates in Hilbert-Schmidt. We complete their work by providing the
rates and optimality for convergence in mean square.

Our initial philosophy is close to the article by Crambes, Kneip and Sarda
(2009). Like these authors we consider the prediction with random design. We
think that this way seems to be the most justified from a statistical point
of view. The case of a fixed design gives birth to several situations and
different rates (with possible oversmoothing which entails parametric rates
of convergence which are odd in this truly nonparametric model) and does not
necessarily correspond to the statistical reality. The main differences rely
in the fact that our results hold in mean square norm rather than in
probability for a larger class of data and parameter at the expense of more
restricted moment assumptions.

Our methodology is closer to the articles by Hall and Horowitz (2007). They
studied the prediction risk at a fixed design in the model with real outputs
(\ref{scalar-model}) but with specified eigenvalues namely $\lambda _{j}\sim
Cj^{-1-\alpha }$ and parameter spectral decomposition $\left\langle \beta
,e_{j}\right\rangle \sim Cj^{-1-\gamma }$ with $\alpha ,\gamma >0$. The
comparisons may be simpler with these works since we share the approach
through spectral decomposition of operators or Karhunen-Loeve development
for the design $X.$

The problem of weak convergence is considered only in Yao, M\"{u}ller and
Wang (2005) : they provide very useful and practical pointwise confidence
sets which imply estimation of the covariance of the noise. Our result may
allow to consider a larger class of testing issues through delta-methods (we
have in mind testing of hypotheses like $S=S_{0}$ versus $S_{\left( n\right)
}=S_{0}+\eta _{n}v$ where $\eta _{n}\rightarrow 0$ and $v$ belongs to a
well-chosen set in $H$).

The contribution of this article essentially deals with a linear regression
model -the concerns related to the functional outputs concentrate on lower
bounds in optimality results and in proving weak convergence with specific
techniques adapted to functional data. We hope that our methods will
demonstrate that optimal results are possible in a general framework and
that regularity assumptions can often be relaxed thanks to the compensation
(or regularity/inverse problem trade-off) phenomenon mentioned within Remark %
\ref{tradeoff}. The Hilbert space framework is necessary at least in the
section devoted to weak convergence. Generalizations to Banach spaces of
functions could be investigated, for instance in $C\left( \left[ 0,1\right]
\right) ,$ H\"{o}lder or Besov spaces.

Finally we do not investigate in this paper the practical point of view of
this prediction method. It is a work in progress. Many directions can be
considered. The practical choice of $k_{n}$ is crucial. Since we provide the
exact theoretical formula for the optimal projection dimension at (\ref%
{k-opt}) it would be interesting to compare it with the results of a
cross-validation method on a simulated dataset. The covariance structure of the noise is a central and major concern : the
covariance operator appears in the limiting distribution, its trace
determines the optimal choice of the dimension $k_{n}^{\ast}.$ Estimating $%
\Gamma _{\varepsilon}$ turns out to be challenging both from a practical and
applied point of view.

\section{Mathematical derivations}

In the sequel, the generic notation $C$ stands for a constant which does not
depend on $k,$ $n$ or $S.$ All our results are related to the decomposition
given below :%
\begin{equation}
\widehat{S}_{n}=S\Gamma _{n}\Gamma _{n}^{\dag }+U_{n}\Gamma _{n}^{\dag }=S%
\widehat{\Pi }_{k}+\frac{1}{n}\sum_{i=1}^{n}\varepsilon _{i}\otimes \Gamma
_{n}^{\dag }X_{i}.  \label{decomp}
\end{equation}%
It is plain that a bias-variance decomposition is exhibited just above. The
random projection $\widehat{\Pi }_{k}$ is not a satisfactory term and we
intend to remove it and to replace it with its non-random counterpart. When
turning to the predictor, (\ref{decomp}) may be enhanced :%
\begin{align}
& \widehat{S}_{n}\left( X_{n+1}\right) -S\left( X_{n+1}\right)
\label{decomp-pred} \\
& =S\left( \Pi _{k}-I\right) \left( X_{n+1}\right) +S\left[ \widehat{\Pi }%
_{k}-\Pi _{k}\right] \left( X_{n+1}\right) +\frac{1}{n}\sum_{i=1}^{n}%
\varepsilon _{i}\left\langle \Gamma _{n}^{\dag }X_{i},X_{n+1}\right\rangle ,
\notag
\end{align}

\noindent where $\Pi_{k}$ is defined in the same way as we defined $\widehat{%
\Pi}_{k}$ previously, \textit{i.e.} the projection on the $k$ first
eigenvectors of $\Gamma$.

In terms of mean square error, the following easily stems from $\mathbb{E}%
\left( \varepsilon_{i}|X\right) =0$~:%
\begin{align*}
& \mathbb{E}\left\Vert \widehat{S}_{n}\left( X_{n+1}\right) -S\left(
X_{n+1}\right) \right\Vert ^{2} \\
& =\mathbb{E}\left\Vert S\widehat{\Pi}_{k}\left( X_{n+1}\right) -S\left(
X_{n+1}\right) \right\Vert ^{2}+\mathbb{E}\left\Vert \frac{1}{n}\sum
_{i=1}^{n}\varepsilon_{i}\left\langle
\Gamma_{n}^{\dag}X_{i},X_{n+1}\right\rangle \right\Vert ^{2}.
\end{align*}

We prove below that :%
\begin{equation}
\mathbb{E}\left\Vert S\left[ \widehat{\Pi}_{k}-\Pi_{k}\right] \left(
X_{n+1}\right) \right\Vert ^{2}=o\left( \mathbb{E}\left\Vert \frac{1}{n}%
\sum_{i=1}^{n}\varepsilon_{i}\left\langle
\Gamma_{n}^{\dag}X_{i},X_{n+1}\right\rangle \right\Vert ^{2}\right),
\label{interm}
\end{equation}
and that the two terms that actually influence the mean square error are the
first and the third in display (\ref{decomp-pred}). The first term $S\left(
\Pi_{k}-I\right) \left( X_{n+1}\right) $ is the bias term and the third a
variance term (see display (\ref{mse})).

The proofs are split into two parts. In the first, part we provide some
technical lemmas which are collected there to enhance the reading of the
second part devoted to the proof of the main results. In all the sequel, the
sequence $k=k_{n}$ depends on $n$ even if this index is dropped. We assume
that all the assumptions mentioned earlier in the paper hold ; they will be
however recalled when addressing crucial steps. We assume once and for all
that $\left( k\log k\right) ^{2}/n\rightarrow 0$ as announced in Remark \ref%
{conv-k} above. The rate of convergence to $0$ of $\left( k\log k\right)
^{2}/n$ will be tuned when dealing with weak convergence.

\subsection{Preliminary material}

\bigskip

All along the proofs, we will make an intensive use of perturbation theory
for bounded operators. It may be useful to have basic notions about spectral
representation of bounded operators and perturbation theory. We refer to
Kato (1976), Dunford and Schwartz (1988, Chapter VII.3) or to Gohberg,
Goldberg and Kaashoek (1991) for an introduction to functional calculus for
operators related with Riesz integrals. Roughly speaking, several results
mentioned below and throughout the article may be easily understood by
considering the formula of residues for analytic functions on the complex
plane (see Rudin (1987)) and extending it to functions still defined on the
complex plane but with values in the space of operators. The introduction of
Gohberg, Goldberg and Kaashoek (1991, pp. 4-16) is illuminating with respect
to this issue.

Let us denote by $\mathcal{B}_{j}$ the oriented circle of the complex plane
with center $\lambda _{j}$ and radius $\delta _{j}/2$ where $\delta
_{j}=\min \left\{ \lambda _{j}-\lambda _{j+1},\lambda _{j-1}-\lambda
_{j}\right\} =\lambda _{j}-\lambda _{j+1}$, the last equality coming from
the convexity associated to the $\lambda _{j}$'s. Let us define $\mathcal{C}%
_{k}=\bigcup_{j=1}^{k}\mathcal{B}_{j}\ .$The open domain whose boundary is $%
\mathcal{C}_{k}$ is not connected but we can apply the functional calculus
for bounded operators (see Dunford-Schwartz, Section VII.3, Definitions 8
and 9). With this formalism at hand it is easy to prove the following
formulas :%
\begin{align}
\Pi _{k_{n}}& =\frac{1}{2\pi \iota }\int_{\mathcal{C}_{k}}\left( zI-\Gamma
\right) ^{-1}dz, \\
\Gamma ^{\dag }& =\frac{1}{2\pi \iota }\int_{\mathcal{C}_{k}}\frac{1}{z}%
\left( zI-\Gamma \right) ^{-1}dz.  \label{resinv}
\end{align}

The same is true with the random $\Gamma _{n}$, but the contour $\mathcal{C}%
_{k}$ must be replaced by its random counterpart $\widehat{\mathcal{C}}%
_{k}=\bigcup_{j=1}^{k_{n}}\widehat{\mathcal{B}}_{j}$ where each $\widehat{%
\mathcal{B}}_{j}$ is a random ball of the complex plane with center $%
\widehat{\lambda }_{j}$ and for instance a radius $\widehat{\delta }_{j}/2$
with plain notations. Then 
\begin{equation*}
\widehat{\Pi }_{k_{n}}=\frac{1}{2\pi \iota }\int_{\widehat{\mathcal{C}}%
_{k}}\left( zI-\Gamma _{n}\right) ^{-1}dz,\quad \Gamma _{n}^{\dag }=\frac{1}{%
2\pi \iota }\int_{\widehat{\mathcal{C}}_{k}}\frac{1}{z}\left( zI-\Gamma
_{n}\right) ^{-1}dz.
\end{equation*}

This first lemma is based on convex inequalities. In the sequel, much
depends on the bounds derived in this Lemma.

\begin{lemma}
\label{L1}Consider two large enough positive integers $j$ and $k$ such that $%
k>j$. Then 
\begin{gather}
j\lambda _{j}\ \geq \ k\lambda _{k},\quad \lambda _{j}-\lambda _{k}\geq
\left( 1-\frac{j}{k}\right) \lambda _{j},\quad \sum_{j\geq k}\lambda
_{j}\leq \left( k+1\right) \lambda _{k}.  \label{t2} \\
\sum_{j\geq 1,j\neq k}\frac{\lambda _{j}}{\left\vert \lambda _{k}-\lambda
_{j}\right\vert }\leq Ck\log k.  \notag
\end{gather}%
Besides 
\begin{equation*}
\mathbb{E}\sup_{z\in \mathcal{B}_{j}}\left\Vert \left( zI-\Gamma \right)
^{-1/2}\left( \Gamma -\Gamma _{n}\right) \left( zI-\Gamma \right)
^{-1/2}\right\Vert _{\mathcal{L}_{2}}^{2}\leq \frac{C}{n}\left( j\log
j\right) ^{2}.
\end{equation*}
\end{lemma}

The proof of this lemma will be found in Cardot, Mas, Sarda (2007), pp.
339-342.

We introduce the following event :%
\begin{equation*}
\mathcal{A}_{n}=\left\{ \forall j\in \left\{ 1,...,k_{n}\right\} ,\;\;\frac{%
\left\vert \widehat{\lambda }_{j}-\lambda _{j}\right\vert }{\delta _{j}}%
<1/2\right\} .
\end{equation*}%
which decribes the way the estimated eigenvalues concentrate around the
population ones : the higher the index $j$ the closer are the $\widehat{%
\lambda }_{j}$'s to the $\lambda _{j}$'s.

\begin{proposition}
\label{P1}If $\left( k\log k\right) ^{2}/n\rightarrow 0$,%
\begin{equation*}
\mathbb{P}\left( \lim \sup \overline{\mathcal{A}}_{n}\right) =0.
\end{equation*}
\end{proposition}

\textbf{Proof :} We just check that the Borel-Cantelli lemma holds $%
\sum_{n=1}^{+\infty }\mathbb{P}\left( \overline{\mathcal{A}}_{n}\right)
<+\infty $ where%
\begin{align*}
\mathbb{P}\left( \overline{\mathcal{A}}_{n}\right) & =\mathbb{P}\left(
\exists j\in \left\{ 1,...,k_{n}\right\} |\left\vert \widehat{\lambda }%
_{j}-\lambda _{j}\right\vert /\delta _{j}>1/2\right) \\
& \leq \sum_{j=1}^{k}\mathbb{P}\left( \left\vert \widehat{\lambda }%
_{j}-\lambda _{j}\right\vert /\lambda _{j}>\delta _{j}/\left( 2\lambda
_{j}\right) \right) \leq \sum_{j=1}^{k}\mathbb{P}\left( \left\vert \widehat{%
\lambda }_{j}-\lambda _{j}\right\vert /\lambda _{j}>1/2\left( j+1\right)
\right) .
\end{align*}%
Now, applying the asymptotic results proved in Bosq (2000) at page 122-124,
we see that the asymptotic behaviour of $\mathbb{P}\left( \left\vert 
\widehat{\lambda }_{j}-\lambda _{j}\right\vert /\lambda _{j}>\frac{1}{2j}%
\right) $ is the same as 
\begin{equation*}
\mathbb{P}\left( \left\vert \frac{1}{n}\sum_{i=1}^{n}\left\langle
X_{i},e_{j}\right\rangle ^{2}-\lambda _{j}\right\vert >\frac{\lambda _{j}}{%
2\left( j+1\right) }\right) .
\end{equation*}%
We apply Bernstein's exponential inequality -which is possible due to
assumption (\ref{assumpt-bernstein})- to the latter, and we obtain (for the
sake of brevity $j+1$ was replaced by $j$ in the right side of the
probability but this does not change the final result) :%
\begin{equation*}
\mathbb{P}\left( \left\vert \frac{1}{n}\sum_{i=1}^{n}\left\langle
X_{i},e_{j}\right\rangle ^{2}-\lambda _{j}\right\vert >\frac{\lambda _{j}}{2j%
}\right) \leq 2\exp \left( -\frac{n}{j^{2}}\frac{1}{8c+1/\left( 6j\right) }%
\right) \leq 2\exp \left( -C\frac{n}{j^{2}}\right) ,
\end{equation*}%
and then%
\begin{equation*}
\sum_{j=1}^{k}\mathbb{P}\left( \left\vert \widehat{\lambda }_{j}-\lambda
_{j}\right\vert >\frac{\lambda _{j}}{2j}\right) \leq 2k\exp \left( -C\frac{n%
}{k^{2}}\right) .
\end{equation*}%
Now it is plain from $\left( k\log k\right) ^{2}/n\rightarrow 0$ that $k\exp
\left( -C\frac{n}{k^{2}}\right) \leq 1/n^{1+\varepsilon }$ for some $%
\varepsilon >0$ which leads to checking that $\sum_{n}k_{n}\exp \left( -C%
\frac{n}{k_{n}^{2}}\right) <+\infty ,$ and to the statement of Proposition %
\ref{P1} through Borel-Cantelli's Lemma.

\begin{corollary}
\label{coro2} We may write%
\begin{equation*}
\widehat{\Pi }_{k_{n}}=\frac{1}{2\pi \iota }\int_{\mathcal{C}_{k}}\left(
zI-\Gamma _{n}\right) ^{-1}dz,\quad \Gamma _{n}^{\dag }=\frac{1}{2\pi \iota }%
\int_{\mathcal{C}_{k}}\frac{1}{z}\left( zI-\Gamma _{n}\right) ^{-1}dz\quad
a.s.,
\end{equation*}%
where this time the contour is $\mathcal{C}_{k}$ hence no more random.
\end{corollary}

\textbf{Proof :} From Proposition \ref{P1}, it is plain that we may assume
that almost surely $\widehat{\lambda}_{j}\in\mathcal{B}_{j}$ for $%
j\in\left\{ 1,...,k\right\} .$ Then the formulas above easily stem from
perturbation theory (see Kato (1976), Dunford and Schwartz (1988) for
instance).

\subsection{Proofs of the main results}

We start with proving (\ref{interm}) as announced in the foreword of this
section. What we give here is nothing but the term $A_{n}$ in Theorem \ref%
{TH2}.

\begin{proposition}
\label{ks}The following bound holds :%
\begin{equation*}
\mathbb{E}\left\Vert S\left( \widehat{\Pi}_{k}-\Pi_{k}\right) \left(
X_{n+1}\right) \right\Vert ^{2}\leq C\frac{k^{2}\lambda_{k}}{n}\left\Vert
S\right\Vert _{\mathcal{L}_{2}}.
\end{equation*}
\end{proposition}

\textbf{Proof : }We start with noting that%
\begin{align*}
\mathbb{E}\left\Vert S\left( \widehat{\Pi }_{k}-\Pi _{k}\right) \left(
X_{n+1}\right) \right\Vert ^{2}& =\mathbb{E}\left[ \mathrm{tr}\left( \Gamma
\left( \widehat{\Pi }_{k}-\Pi _{k}\right) S^{\ast }S\left( \widehat{\Pi }%
_{k}-\Pi _{k}\right) \right) \right] \\
& =\mathbb{E}\left\Vert S\left( \widehat{\Pi }_{k}-\Pi _{k}\right) \Gamma
^{1/2}\right\Vert _{\mathcal{L}_{2}}^{2} \\
& =\sum_{j=1}^{+\infty }\sum_{\ell =1}^{+\infty }\left\langle S\left( 
\widehat{\Pi }_{k}-\Pi _{k}\right) \Gamma ^{1/2}\left( e_{j}\right) ,e_{\ell
}\right\rangle ^{2}.
\end{align*}

By Corollary \ref{coro2}, we have 
\begin{equation}
\widehat{\Pi }_{k}-\Pi _{k}=\frac{1}{2\pi \iota }\sum_{m=1}^{k}\int_{%
\mathcal{B}_{m}}\left\{ \left( zI-\Gamma _{n}\right) ^{-1}-\left( zI-\Gamma
\right) ^{-1}\right\} dz=\sum_{m=1}^{k}T_{m,n},  \label{sw}
\end{equation}%
where $T_{m,n}=\frac{1}{2\pi \iota }\int_{\mathcal{B}_{m}}\left( zI-\Gamma
_{n}\right) ^{-1}\left( \Gamma -\Gamma _{n}\right) \left( zI-\Gamma \right)
^{-1}dz$.

To go ahead now, we ask the reader to accept momentaneously that for all $%
m\leq k$, the asymptotic behaviour of $T_{m,n}$ is the same as 
\begin{equation*}
T_{m,n}^{\ast }=\frac{1}{2\pi \iota }\int_{\mathcal{B}_{m}}\left( zI-\Gamma
\right) ^{-1}\left( \Gamma -\Gamma _{n}\right) \left( zI-\Gamma \right)
^{-1}dz,
\end{equation*}%
where the random $\left( zI-\Gamma _{n}\right) ^{-1}$ was replaced by the
non-random $\left( zI-\Gamma \right) ^{-1}$ and that studying $\widehat{\Pi }%
_{k}-\Pi _{k}$ comes down to studying 
\begin{equation*}
\frac{1}{2\pi \iota }\sum_{m=1}^{k}\int_{\mathcal{B}_{m}}\left( zI-\Gamma
\right) ^{-1}\left( \Gamma -\Gamma _{n}\right) \left( zI-\Gamma \right)
^{-1}dz.
\end{equation*}%
The proof that this switch is allowed is postponed to Lemma \ref{switch}. We
go on with%
\begin{align*}
& \left\langle S\left( \widehat{\Pi }_{k}-\Pi _{k}\right) \Gamma
^{1/2}\left( e_{j}\right) ,e_{\ell }\right\rangle =\frac{1}{2\pi \iota }%
\sum_{m=1}^{k}\int_{\mathcal{B}_{m}}\left\langle \left( zI-\Gamma \right)
^{-1}\left( \Gamma -\Gamma _{n}\right) \left( zI-\Gamma \right) ^{-1}\Gamma
^{1/2}\left( e_{j}\right) ,S^{\ast }e_{\ell }\right\rangle dz \\
& =\frac{\sqrt{\lambda _{j}}}{2\pi \iota }\sum_{m=1}^{k}\int_{\mathcal{B}%
_{m}}\left\langle \left( zI-\Gamma \right) ^{-1}\left( \Gamma -\Gamma
_{n}\right) \left( e_{j}\right) ,S^{\ast }e_{\ell }\right\rangle \frac{dz}{%
z-\lambda _{j}},
\end{align*}%
where $S^{\ast }$ is the adjoint operator of $S$. We obtain 
\begin{align*}
& \int_{\mathcal{B}_{m}}\left\langle \left( zI-\Gamma \right) ^{-1}\left(
\Gamma -\Gamma _{n}\right) \left( e_{j}\right) ,S^{\ast }e_{\ell
}\right\rangle \frac{dz}{z-\lambda _{j}} \\
& =\int_{\mathcal{B}_{m}}\sum_{j^{\prime }=1}^{+\infty }\left\langle \left(
zI-\Gamma \right) ^{-1}\left( \Gamma -\Gamma _{n}\right) \left( e_{j}\right)
,e_{j^{\prime }}\right\rangle \left\langle S^{\ast }e_{\ell },e_{j^{\prime
}}\right\rangle \frac{dz}{z-\lambda _{j}} \\
& =\int_{\mathcal{B}_{m}}\sum_{j^{\prime }=1}^{+\infty }\left\langle \left(
\Gamma -\Gamma _{n}\right) \left( e_{j}\right) ,e_{j^{\prime }}\right\rangle
\left\langle S^{\ast }e_{\ell },e_{j^{\prime }}\right\rangle \frac{dz}{%
\left( z-\lambda _{j}\right) \left( z-\lambda _{j^{\prime }}\right) }.
\end{align*}%
We deduce that%
\begin{equation*}
\left\langle S\left( \widehat{\Pi }_{k}-\Pi _{k}\right) \Gamma ^{1/2}\left(
e_{j}\right) ,e_{\ell }\right\rangle =\frac{\sqrt{\lambda _{j}}}{2\pi \iota }%
\sum_{j^{\prime }=1}^{+\infty }\left\langle \left( \Gamma -\Gamma
_{n}\right) \left( e_{j}\right) ,e_{j^{\prime }}\right\rangle \left\langle
S^{\ast }e_{\ell },e_{j^{\prime }}\right\rangle \sum_{m=1}^{k}\int_{\mathcal{%
B}_{m}}\frac{dz}{\left( z-\lambda _{j}\right) \left( z-\lambda _{j^{\prime
}}\right) },
\end{equation*}%
then%
\begin{align*}
& \left\langle S\left( \widehat{\Pi }_{k}-\Pi _{k}\right) \Gamma
^{1/2}\left( e_{j}\right) ,e_{\ell }\right\rangle \\
& =\frac{\sqrt{\lambda _{j}}}{2\pi \iota }\sum_{j^{\prime
}=1}^{k}\left\langle \left( \Gamma -\Gamma _{n}\right) \left( e_{j}\right)
,e_{j^{\prime }}\right\rangle \left\langle S^{\ast }e_{\ell },e_{j^{\prime
}}\right\rangle \sum_{m=1}^{k}\int_{\mathcal{B}_{m}}\frac{dz}{\left(
z-\lambda _{j}\right) \left( z-\lambda _{j^{\prime }}\right) } \\
& +\frac{\sqrt{\lambda _{j}}}{2\pi \iota }\sum_{j^{\prime }=k+1}^{+\infty
}\left\langle \left( \Gamma -\Gamma _{n}\right) \left( e_{j}\right)
,e_{j^{\prime }}\right\rangle \left\langle S^{\ast }e_{\ell },e_{j^{\prime
}}\right\rangle \sum_{m=1}^{k}\int_{\mathcal{B}_{m}}\frac{dz}{\left(
z-\lambda _{j}\right) \left( z-\lambda _{j^{\prime }}\right) },
\end{align*}%
where%
\begin{equation*}
\sum_{m=1}^{k}\int_{\mathcal{B}_{m}}\frac{dz}{\left( z-\lambda _{j}\right)
\left( z-\lambda _{j^{\prime }}\right) }=\left\{ 
\begin{array}{l}
0\text{ if $j,j^{\prime }>m$}, \\ 
\left( \lambda _{j}-\lambda _{j^{\prime }}\right) ^{-1}\text{ if $j^{\prime
}>m,j\leq m$}, \\ 
\left( \lambda _{j^{\prime }}-\lambda _{j}\right) ^{-1}\text{ if $j^{\prime
}\leq m,j>m$}, \\ 
1-1=0\text{ if $j,j^{\prime }\leq m$}.%
\end{array}%
\right.
\end{equation*}%
Then%
\begin{align*}
& \sum_{j=1}^{+\infty }\left\langle S\left( \widehat{\Pi }_{k}-\Pi
_{k}\right) \Gamma ^{1/2}\left( e_{j}\right) ,e_{\ell }\right\rangle
^{2}=\sum_{j=1}^{k}\left[ \frac{\sqrt{\lambda _{j}}}{2\pi \iota }%
\sum_{j^{\prime }=k+1}^{+\infty }\frac{\left\langle \left( \Gamma -\Gamma
_{n}\right) \left( e_{j}\right) ,e_{j^{\prime }}\right\rangle }{\left(
\lambda _{j}-\lambda _{j^{\prime }}\right) }\left\langle S^{\ast }e_{\ell
},e_{j^{\prime }}\right\rangle \right] ^{2} \\
& +\sum_{j=k+1}^{+\infty }\left[ \frac{\sqrt{\lambda _{j}}}{2\pi \iota }%
\sum_{j^{\prime }=1}^{k}\frac{\left\langle \left( \Gamma -\Gamma _{n}\right)
\left( e_{j}\right) ,e_{j^{\prime }}\right\rangle }{\left( \lambda
_{j^{\prime }}-\lambda _{j}\right) }\left\langle S^{\ast }e_{\ell
},e_{j^{\prime }}\right\rangle \right] ^{2}=A+B,
\end{align*}%
where%
\begin{align*}
A& =\frac{1}{4\pi ^{2}}\sum_{j=1}^{k}\lambda _{j}\left[ \sum_{j^{\prime
}=k+1}^{+\infty }\frac{\left\langle \left( \Gamma -\Gamma _{n}\right) \left(
e_{j}\right) ,e_{j^{\prime }}\right\rangle }{\left( \lambda _{j}-\lambda
_{j^{\prime }}\right) }\left\langle S^{\ast }e_{\ell },e_{j^{\prime
}}\right\rangle \right] ^{2}, \\
B& =\frac{1}{4\pi ^{2}}\sum_{j=k+1}^{+\infty }\lambda _{j}\left[
\sum_{j^{\prime }=1}^{k}\frac{\left\langle \left( \Gamma -\Gamma _{n}\right)
\left( e_{j}\right) ,e_{j^{\prime }}\right\rangle }{\left( \lambda
_{j^{\prime }}-\lambda _{j}\right) }\left\langle S^{\ast }e_{\ell
},e_{j^{\prime }}\right\rangle \right] ^{2}.
\end{align*}%
We first compute $\mathbb{E}A$. To that aim we focus on%
\begin{eqnarray*}
&&\mathbb{E}\left[ \sum_{j^{\prime }=k+1}^{+\infty }\frac{\left\langle
\left( \Gamma -\Gamma _{n}\right) \left( e_{j}\right) ,e_{j^{\prime
}}\right\rangle }{\left( \lambda _{j}-\lambda _{j^{\prime }}\right) }%
\left\langle S^{\ast }e_{\ell },e_{j^{\prime }}\right\rangle \right]
^{2}=\sum_{j^{\prime }=k+1}^{+\infty }\frac{\mathbb{E}\left\langle \left(
\Gamma _{n}-\Gamma \right) \left( e_{j}\right) ,e_{j^{\prime }}\right\rangle
^{2}}{\left( \lambda _{j}-\lambda _{j^{\prime }}\right) ^{2}}\left\langle
S^{\ast }e_{\ell },e_{j^{\prime }}\right\rangle ^{2} \\
&&+\sum_{{\substack{j',j''=k+1 \\ j' \neq j''}}}^{+\infty }\mathbb{E}%
\left\langle \left( \Gamma _{n}-\Gamma \right) \left( e_{j}\right)
,e_{j^{\prime }}\right\rangle \left\langle \left( \Gamma -\Gamma _{n}\right)
\left( e_{j}\right) ,e_{j^{\prime \prime }}\right\rangle \frac{\left\langle
S^{\ast }e_{\ell },e_{j^{\prime }}\right\rangle \left\langle S^{\ast
}e_{\ell },e_{j^{\prime \prime }}\right\rangle }{\left( \lambda _{j}-\lambda
_{j^{\prime }}\right) \left( \lambda _{j}-\lambda _{j^{\prime \prime
}}\right) } \\
&=&\frac{1}{n}\sum_{j^{\prime }=k+1}^{+\infty }c_{j,j^{\prime }}\frac{%
\lambda _{j}\lambda _{j^{\prime }}}{\left( \lambda _{j}-\lambda _{j^{\prime
}}\right) ^{2}}\left\langle S^{\ast }e_{\ell },e_{j^{\prime }}\right\rangle
^{2} \\
&&+\frac{1}{n}\sum_{{\substack{j',j''=k+1 \\ j' \neq j''}}}^{+\infty }%
\mathbb{E}\left[ \left\langle X,e_{j}\right\rangle ^{2}\left\langle
X,e_{j^{\prime }}\right\rangle \left\langle X,e_{j^{\prime \prime
}}\right\rangle \right] \frac{\left\langle S^{\ast }e_{\ell },e_{j^{\prime
}}\right\rangle \left\langle S^{\ast }e_{\ell },e_{j^{\prime \prime
}}\right\rangle }{\left( \lambda _{j}-\lambda _{j^{\prime }}\right) \left(
\lambda _{j}-\lambda _{j^{\prime \prime }}\right) }
\end{eqnarray*}

Then 
\begin{eqnarray*}
&&\mathbb{E}\left[ \sum_{j^{\prime }=k+1}^{+\infty }\frac{\left\langle
\left( \Gamma -\Gamma _{n}\right) \left( e_{j}\right) ,e_{j^{\prime
}}\right\rangle }{\left( \lambda _{j}-\lambda _{j^{\prime }}\right) }%
\left\langle S^{\ast }e_{\ell },e_{j^{\prime }}\right\rangle \right] ^{2} \\
&\leq &C_{1}\frac{\lambda _{j}}{n}\sum_{j^{\prime }=k+1}^{+\infty }\frac{%
\lambda _{j^{\prime }}}{\left( \lambda _{j}-\lambda _{j^{\prime }}\right)
^{2}}\left\langle S^{\ast }e_{\ell },e_{j^{\prime }}\right\rangle ^{2}+C_{2}%
\frac{\lambda _{j}}{n}\sum_{{\substack{j',j''=k+1 \\ j' \neq j''}}}^{+\infty
}\sqrt{\lambda _{j^{\prime }}\lambda _{j^{\prime \prime }}}\frac{%
\left\langle S^{\ast }e_{\ell },e_{j^{\prime }}\right\rangle \left\langle
S^{\ast }e_{\ell },e_{j^{\prime \prime }}\right\rangle }{\left( \lambda
_{j}-\lambda _{j^{\prime }}\right) \left( \lambda _{j}-\lambda _{j^{\prime
\prime }}\right) } \\
&\leq &C\frac{\lambda _{j}}{n}\left( \sum_{j^{\prime }=k+1}^{+\infty }\frac{%
\sqrt{\lambda _{j^{\prime }}}}{\left( \lambda _{j}-\lambda _{j^{\prime
}}\right) }\left\langle S^{\ast }e_{\ell },e_{j^{\prime }}\right\rangle
\right) ^{2}.
\end{eqnarray*}%
We could prove exactly in the same way that%
\begin{equation}
\mathbb{E}\left[ \sum_{j^{\prime }=1}^{k}\frac{\left\langle \left( \Gamma
_{n}-\Gamma \right) \left( e_{j}\right) ,e_{j^{\prime }}\right\rangle }{%
\left( \lambda _{j^{\prime }}-\lambda _{j}\right) }\left\langle S^{\ast
}e_{\ell },e_{j^{\prime }}\right\rangle \right] ^{2}\leq C^{\prime }\frac{%
\lambda _{j}}{n}\left( \sum_{j^{\prime }=1}^{k}\frac{\sqrt{\lambda
_{j^{\prime }}}}{\left( \lambda _{j^{\prime }}-\lambda _{j}\right) }%
\left\langle S^{\ast }e_{\ell },e_{j^{\prime }}\right\rangle \right) ^{2}.
\label{L2}
\end{equation}%
We turn back to 
\begin{align*}
& \left\vert \sum_{j^{\prime }=k+1}^{+\infty }\frac{\sqrt{\lambda
_{j^{\prime }}}}{\left( \lambda _{j}-\lambda _{j^{\prime }}\right) }%
\left\langle S^{\ast }e_{\ell },e_{j^{\prime }}\right\rangle \right\vert
\leq \sum_{j^{\prime }=k+1}^{+\infty }\frac{\sqrt{\lambda _{j^{\prime }}}}{%
\left( \lambda _{j}-\lambda _{j^{\prime }}\right) }\left\vert \left\langle
S^{\ast }e_{\ell },e_{j^{\prime }}\right\rangle \right\vert \\
& =\sum_{j^{\prime }=k+1}^{2k}\frac{\sqrt{\lambda _{j^{\prime }}}}{\left(
\lambda _{j}-\lambda _{j^{\prime }}\right) }\left\vert \left\langle S^{\ast
}e_{\ell },e_{j^{\prime }}\right\rangle \right\vert +\sum_{j^{\prime
}=2k+1}^{+\infty }\frac{\sqrt{\lambda _{j^{\prime }}}}{\left( \lambda
_{j}-\lambda _{j^{\prime }}\right) }\left\vert \left\langle S^{\ast }e_{\ell
},e_{j^{\prime }}\right\rangle \right\vert \\
& \leq \frac{\sqrt{\lambda _{k+1}}}{\left( \lambda _{j}-\lambda
_{k+1}\right) }\sum_{j^{\prime }=k+1}^{2k}\left\vert \left\langle S^{\ast
}e_{\ell },e_{j^{\prime }}\right\rangle \right\vert +\frac{2}{\lambda _{j}}%
\sum_{j^{\prime }=2k+1}^{+\infty }\sqrt{\lambda _{j^{\prime }}}\left\vert
\left\langle S^{\ast }e_{\ell },e_{j^{\prime }}\right\rangle \right\vert ,
\end{align*}%
hence%
\begin{eqnarray}
\mathbb{E}A &\leq &\frac{C}{n}\sum_{j=1}^{k}\lambda _{j}^{2}\left[ \frac{%
\lambda _{k+1}}{\left( \lambda _{j}-\lambda _{k+1}\right) ^{2}}\left(
\sum_{j^{\prime }=k+1}^{2k}\left\vert \left\langle S^{\ast }e_{\ell
},e_{j^{\prime }}\right\rangle \right\vert \right) ^{2}\right]  \label{ring}
\\
&&+.\frac{Ck}{n}\left( \sum_{j^{\prime }=2k+1}^{+\infty }\sqrt{\lambda
_{j^{\prime }}}\left\vert \left\langle S^{\ast }e_{\ell },e_{j^{\prime
}}\right\rangle \right\vert \right) ^{2}  \notag
\end{eqnarray}%
The term below is bounded by :%
\begin{equation*}
\frac{Ck}{n}\left( \sum_{j^{\prime }=2k+1}^{+\infty }\lambda _{j^{\prime
}}\sum_{j^{\prime }=2k+1}^{+\infty }\left\vert \left\langle S^{\ast }e_{\ell
},e_{j^{\prime }}\right\rangle \right\vert ^{2}\right) \leq \frac{Ck^{2}}{n}%
\lambda _{k}\sum_{j^{\prime }=2k+1}^{+\infty }\left\vert \left\langle
S^{\ast }e_{\ell },e_{j^{\prime }}\right\rangle \right\vert ^{2}
\end{equation*}%
because $\sum_{j^{\prime }=2k+1}^{+\infty }\lambda _{j^{\prime }}\leq \left(
2k+1\right) \lambda _{2k+1}\leq k\lambda _{k}$ by Lemma \ref{L1}. We focus
on the term on line (\ref{ring}):%
\begin{align*}
& \sum_{j=1}^{k}\lambda _{j}^{2}\left[ \frac{\lambda _{k+1}}{\left( \lambda
_{j}-\lambda _{k+1}\right) ^{2}}\left( \sum_{j^{\prime }=k+1}^{2k}\left\vert
\left\langle S^{\ast }e_{\ell },e_{j^{\prime }}\right\rangle \right\vert
\right) ^{2}\right] \leq \lambda _{k+1}\sum_{j=1}^{k}\left[ \left( \frac{k+1%
}{k+1-j}\right) ^{2}\left( \sum_{j^{\prime }=k+1}^{2k}\left\vert
\left\langle S^{\ast }e_{\ell },e_{j^{\prime }}\right\rangle \right\vert
\right) ^{2}\right] \\
& \leq k\left( \sum_{j^{\prime }=k+1}^{2k}\left\vert \left\langle S^{\ast
}e_{\ell },e_{j^{\prime }}\right\rangle \right\vert ^{2}\right) \left(
k+1\right) ^{2}\lambda _{k+1}\sum_{j=1}^{k}\frac{1}{j^{2}}\leq C\left(
\sum_{j^{\prime }=k+1}^{2k}\left\vert \left\langle S^{\ast }e_{\ell
},e_{j^{\prime }}\right\rangle \right\vert ^{2}\right) k^{2}\lambda _{k+1},
\end{align*}%
hence $\mathbb{E}A\leq \frac{C}{n}\left( \sum_{j^{\prime }=k+1}^{+\infty
}\left\vert \left\langle S^{\ast }e_{\ell },e_{j^{\prime }}\right\rangle
\right\vert ^{2}\right) k^{2}\lambda _{k}.$ We turn to proving a similar
bound for $B$. The method is given because it is significantly distinct. We
start from (\ref{L2}) and we denote $\lfloor x\rfloor $ the largest integer
smaller than $x$ :%
\begin{align*}
& \frac{\lambda _{j}}{n}\left( \sum_{j^{\prime }=1}^{k}\frac{\sqrt{\lambda
_{j^{\prime }}}}{\left( \lambda _{j^{\prime }}-\lambda _{j}\right) }%
\left\langle S^{\ast }e_{\ell },e_{j^{\prime }}\right\rangle \right) ^{2} \\
& \leq \frac{\lambda _{j}}{n}\left[ \left( \sum_{j^{\prime
}=1}^{\left\lfloor k/2\right\rfloor }\frac{\sqrt{\lambda _{j^{\prime }}}}{%
\left( \lambda _{j^{\prime }}-\lambda _{j}\right) }\left\vert \left\langle
S^{\ast }e_{\ell },e_{j^{\prime }}\right\rangle \right\vert \right)
^{2}+\left( \sum_{j^{\prime }=\left\lfloor k/2\right\rfloor }^{k}\frac{\sqrt{%
\lambda _{j^{\prime }}}}{\left( \lambda _{j^{\prime }}-\lambda _{j}\right) }%
\left\vert \left\langle S^{\ast }e_{\ell },e_{j^{\prime }}\right\rangle
\right\vert \right) ^{2}\right] \\
& \leq C\frac{\lambda _{j}}{n}\left[ \left( \sum_{j^{\prime
}=1}^{\left\lfloor k/2\right\rfloor }\frac{1}{\sqrt{\lambda _{j^{\prime }}}}%
\left\vert \left\langle S^{\ast }e_{\ell },e_{j^{\prime }}\right\rangle
\right\vert \right) ^{2}+\frac{1}{\lambda _{k}-\lambda _{j}}\frac{j}{j-k}%
k\sum_{j^{\prime }=\left\lfloor k/2\right\rfloor }^{k}\left\langle S^{\ast
}e_{\ell },e_{j^{\prime }}\right\rangle ^{2}\right] \\
& \leq C\frac{\lambda _{j}k}{n\lambda _{k}}\sum_{j^{\prime
}=1}^{\left\lfloor k/2\right\rfloor }\left\langle S^{\ast }e_{\ell
},e_{j^{\prime }}\right\rangle ^{2}+\frac{1}{n}\frac{\lambda _{j}}{\lambda
_{k}-\lambda _{j}}\frac{j}{j-k}k\sum_{j^{\prime }=\left\lfloor
k/2\right\rfloor }^{k}\left\langle S^{\ast }e_{\ell },e_{j^{\prime
}}\right\rangle ^{2} \\
& \leq C\frac{k}{n}\sum_{j^{\prime }=1}^{k}\left\langle S^{\ast }e_{\ell
},e_{j^{\prime }}\right\rangle ^{2}+\frac{k}{n}\left( \frac{j}{j-k}\right)
^{2}\sum_{j^{\prime }=\left\lfloor k/2\right\rfloor }^{k}\left\langle
S^{\ast }e_{\ell },e_{j^{\prime }}\right\rangle ^{2}.
\end{align*}%
>From the definition of $B$, we get finally%
\begin{equation*}
\mathbb{E}B\leq C\frac{k}{n}\left( \sum_{j^{\prime }=1}^{k}\left\vert
\left\langle S^{\ast }e_{\ell },e_{j^{\prime }}\right\rangle \right\vert
^{2}\right) \sum_{j=k+1}^{+\infty }\lambda _{j}+\left( \sum_{j^{\prime
}=\left\lfloor k/2\right\rfloor }^{k}\left\langle S^{\ast }e_{\ell
},e_{j^{\prime }}\right\rangle ^{2}\right) \frac{k}{n}\sum_{j=k+1}^{+\infty
}\lambda _{j}\left( \frac{j}{j-k}\right) ^{2}.
\end{equation*}%
It is plain that, for sufficiently large $k$, $\sum_{j^{\prime
}=\left\lfloor k/2\right\rfloor }^{k}\left\langle S^{\ast }e_{\ell
},e_{j^{\prime }}\right\rangle ^{2}\leq C/k$ (otherwise $\sum_{j^{\prime
}}\left\langle S^{\ast }e_{\ell },e_{j^{\prime }}\right\rangle ^{2}$ cannot
converge), whence%
\begin{align*}
& \left( \sum_{j^{\prime }=\left\lfloor k/2\right\rfloor }^{k}\left\langle
S^{\ast }e_{\ell },e_{j^{\prime }}\right\rangle ^{2}\right) \frac{k}{n}%
\sum_{j=k+1}^{+\infty }\lambda _{j}\left( \frac{j}{j-k}\right) ^{2}\leq 
\frac{C}{n}\left[ \sum_{j=k+1}^{2k}\lambda _{j}\left( \frac{j}{j-k}\right)
^{2}+\sum_{j=2k}^{+\infty }\lambda _{j}\left( \frac{j}{j-k}\right) ^{2}%
\right] \\
& \leq \frac{C}{n}\left[ \sum_{j=k+1}^{2k}\lambda _{j}\left( \frac{j}{j-k}%
\right) ^{2}+4\sum_{j=2k}^{+\infty }\lambda _{j}\right] .
\end{align*}%
Denoting $\varkappa _{k}=\sup_{k+1\leq j\leq 2k}\left( j\log j\lambda
_{j}\right) $ we get at last :%
\begin{align*}
\sum_{j=k+1}^{2k}\lambda _{j}\left( \frac{j}{j-k}\right) ^{2}& \leq
\sup_{k+1\leq j\leq 2k}\left( j\log j\lambda _{j}\right) \frac{1}{\log k}%
\sum_{j=k+1}^{2k}\frac{j}{j-k} \\
& \leq \varkappa _{k}\frac{1}{\log k}\sum_{j=1}^{k}\frac{k+j}{j}\leq
Ck\varkappa _{k},
\end{align*}%
and $\mathbb{E}B\leq C\frac{k}{n}\varkappa _{k}\left( \sum_{j^{\prime
}=1}^{k}\left\vert \left\langle S^{\ast }e_{\ell },e_{j^{\prime
}}\right\rangle \right\vert ^{2}\right) ,$with $\varkappa _{k}\rightarrow 0$%
. Finally :%
\begin{equation*}
\sum_{j=1}^{+\infty }\sum_{\ell =1}^{+\infty }\left\langle S\left( \widehat{%
\Pi }_{k}-\Pi _{k}\right) \Gamma ^{1/2}\left( e_{j}\right) ,e_{\ell
}\right\rangle ^{2}\leq C\frac{k}{n}\varkappa _{k}\sum_{j=1}^{+\infty
}\sum_{\ell =1}^{+\infty }\left\vert \left\langle S^{\ast }e_{\ell
},e_{j}\right\rangle \right\vert ^{2}.
\end{equation*}

This last bound almost concludes the rather long proof of Proposition \ref%
{ks}. It remains to ensure that switching $T_{m,n}^{\ast }$ and $T_{m,n}$ as
announced just below display (\ref{sw}) is possible.

\begin{lemma}
\label{switch}We have%
\begin{equation*}
\mathbb{E}\sum_{j=1}^{+\infty}\sum_{\ell=1}^{+\infty}\left\langle S\left( 
\widehat{\Pi}_{k}-\Pi_{k}\right) \Gamma^{1/2}\left( e_{j}\right)
,e_{\ell}\right\rangle ^{2}\sim\mathbb{E}\sum_{j=1}^{+\infty}\sum_{%
\ell=1}^{+\infty }\sum_{m=1}^{k}\left\langle
ST_{m,n}^{\ast}\Gamma^{1/2}\left( e_{j}\right) ,e_{\ell}\right\rangle ^{2}.
\end{equation*}
In other words, switching $T_{m,n}^{\ast}$ and $T_{m,n}$ is possible in
display (\ref{sw}).
\end{lemma}

The proof of this Lemma is close to the control of second order term at page
351-352 of Cardot, Mas and Sarda (2007) and we will give a sketch of it. We
start from :%
\begin{align*}
T_{m,n}& =\frac{1}{2\pi \iota }\int_{\mathcal{B}_{m}}\left( zI-\Gamma
_{n}\right) ^{-1}\left( \Gamma -\Gamma _{n}\right) \left( zI-\Gamma \right)
^{-1}dz \\
& =\frac{1}{2\pi \iota }\int_{\mathcal{B}_{m}}\left( zI-\Gamma \right)
^{-1/2}R_{n}\left( z\right) \left( zI-\Gamma \right) ^{-1/2}\left( \Gamma
-\Gamma _{n}\right) \left( zI-\Gamma \right) ^{-1}dz,
\end{align*}%
with $R_{n}\left( z\right) =\left( zI-\Gamma \right) ^{1/2}\left( zI-\Gamma
_{n}\right) ^{-1}\left( zI-\Gamma \right) ^{1/2}$. Besides, as can be seen
from Lemma 4 in\ Cardot, Mas and Sarda (2007)

\begin{equation*}
\left[ I+\left( zI-\Gamma \right) ^{-1/2}\left( \Gamma -\Gamma _{n}\right)
\left( zI-\Gamma \right) ^{-1/2}\right] R_{n}\left( z\right) =I.
\end{equation*}%
Denoting $S_{n}\left( z\right) =\left( zI-\Gamma \right) ^{-1/2}\left(
\Gamma -\Gamma _{n}\right) \left( zI-\Gamma \right) ^{-1/2}$, it is plain
that when $\left\Vert S_{n}\left( z\right) \right\Vert \leq 1$ for all $z\in 
\mathcal{C}_{k}$,%
\begin{equation*}
R_{n}\left( z\right) =I+\sum_{m=1}^{+\infty }\left( -1\right)
^{m}S_{n}^{m}\left( z\right) :=I+R_{n}^{0}\left( z\right) ,
\end{equation*}%
with $\left\Vert R_{n}^{0}\left( z\right) \right\Vert _{\infty }\leq
C\left\Vert S_{n}\left( z\right) \right\Vert _{\infty }$ for all $z\in 
\mathcal{C}_{k}.$ Turning back to our initial equation we get, conditionally
to $\left\Vert S_{n}\left( z\right) \right\Vert \leq 1$ for all $z\in 
\mathcal{C}_{k}$ :%
\begin{equation*}
T_{m,n}-T_{m,n}^{\ast }=\frac{1}{2\pi \iota }\int_{\mathcal{B}_{m}}\left(
zI-\Gamma \right) ^{-1/2}R_{n}^{0}\left( z\right) \left( zI-\Gamma \right)
^{-1/2}\left( \Gamma -\Gamma _{n}\right) \left( zI-\Gamma \right) ^{-1}dz,
\end{equation*}%
and we confine to considering only the first term in the devlopment of $%
R_{n}^{0}\left( z\right) $ which writes $\left( 2\pi \iota \right)
^{-1}\int_{\mathcal{B}_{m}}\left( zI-\Gamma \right) ^{-1/2}S_{n}^{2}\left(
z\right) \left( zI-\Gamma \right) ^{-1/2}dz$.

Now split $S\left( \widehat{\Pi }_{k}-\Pi _{k}\right) \Gamma ^{1/2}=S\left( 
\widehat{\Pi }_{k}-\Pi _{k}\right) \Gamma ^{1/2}1\mathrm{I}_{\mathcal{J}%
}+S\left( \widehat{\Pi }_{k}-\Pi _{k}\right) \Gamma ^{1/2}1\mathrm{I}_{%
\overline{\mathcal{J}}}$ where 
\begin{equation*}
\mathcal{J}=\left\{ \sup_{z\in \mathcal{C}_{k}}\left\Vert
\left( zI-\Gamma \right) ^{-1/2}\left( \Gamma -\Gamma _{n}\right) \left(
zI-\Gamma \right) ^{-1/2}\right\Vert _{\mathcal{L}_{2}}^{2}<\tau
_{n}k_{n}/n\right\}
\end{equation*}

and $\tau _{n}$ will be tuned later. We have :%
\begin{equation}
\mathbb{E}\left\Vert S\left( \widehat{\Pi }_{k}-\Pi _{k}\right) \Gamma
^{1/2}1\mathrm{I}_{\overline{\mathcal{J}}}\right\Vert _{\mathcal{L}%
_{2}}^{2}\leq 4\left\Vert S\Gamma ^{1/2}\right\Vert _{\mathcal{L}_{2}}^{2}%
\mathbb{P}\left( \overline{\mathcal{J}}\right),  \label{rem}
\end{equation}%
and 
\begin{eqnarray*}
&&\left\Vert S\left[ \left( \widehat{\Pi }_{k}-\Pi _{k}\right)
-\sum_{m=1}^{k}T_{m,n}^{\ast }\right] \Gamma ^{1/2}1\mathrm{I}_{\mathcal{J}%
}\right\Vert _{\mathcal{L}_{2}} \\
&\leq &\left\Vert S\left[ \sum_{m=1}^{k}\left( 2\pi \iota \right) ^{-1}\int_{%
\mathcal{B}_{m}}\left( zI-\Gamma \right) ^{-1/2}S_{n}^{2}\left( z\right)
\left( zI-\Gamma \right) ^{-1/2}dz\right] \Gamma ^{1/2}1\mathrm{I}_{\mathcal{%
J}}\right\Vert _{\mathcal{L}_{2}} \\
&\leq &\left( 2\pi \right) ^{-1}\frac{\tau _{n}^{2}k_{n}^{2}}{n^{2}}%
\sum_{m=1}^{k}\delta _{m}\sup_{z\in \mathcal{B}_{m}}\left\{ \left\Vert
\left( zI-\Gamma \right) ^{-1/2}\Gamma ^{1/2}\right\Vert _{\infty
}\left\Vert S\left( zI-\Gamma \right) ^{-1/2}\right\Vert _{\infty }\right\}
\\
&\leq &\left( 2\pi \right) ^{-1}\left\Vert S\right\Vert _{\infty }\frac{\tau
_{n}^{2}k_{n}^{2}}{n^{2}}\sum_{m=1}^{k}\sqrt{\delta _{m}m}.
\end{eqnarray*}

Now from $\sum_{m=1}^{+\infty }m\delta _{m}<+\infty $ we get $\sqrt{\delta
_{m}m}\leq c/\sqrt{m\log m}$ hence $\frac{\tau _{n}^{2}k_{n}^{2}}{n^{2}}%
\sum_{m=1}^{k}\sqrt{\delta _{m}m}=o\left( \sqrt{k_{n}/n}\right) $ whenever $%
k_{n}^{4}\tau _{n}^{4}/n^{3}\rightarrow 0.$

The last step consists in controlling the right hand side of (\ref{rem}). In
Cardot, Mas and Sarda (2007) this is done by classical Markov moment
assumptions under the condition that $k_{n}^{5}\log ^{4}n/n$ tends to zero.
Here, Bernstein's exponential inequality yields a tighter bound and ensures
that $\mathbb{P}\left( \overline{\mathcal{J}}\right) =o\left( k_{n}/n\right) 
$ when $k_{n}^{2}\log ^{2}k_{n}/n$ tends to zero. The method of proof is
close in spirit though slightly more intricate than Proposition \ref{P1}.

\begin{proposition}
\label{variance}Let $T_{n}=\frac{1}{n}\sum_{i=1}^{n}\varepsilon
_{i}\left\langle \Gamma _{n}^{\dag }X_{i},X_{n+1}\right\rangle $, then%
\begin{equation*}
\mathbb{E}\left\Vert T_{n}\right\Vert ^{2}=\frac{\sigma _{\varepsilon }^{2}}{%
n}k+\frac{\mathrm{tr}\left[ \Gamma \mathbb{E}\left( \Gamma _{n}^{\dag
}-\Gamma ^{\dag }\right) \right] }{n}.
\end{equation*}
\end{proposition}

\begin{remark}
We see that the right hand side in the display above matches the
decomposition in (\ref{mse}) and $\mathrm{tr}\left[ \Gamma\mathbb{E}\left(
\Gamma _{n}^{\dag}-\Gamma^{\dag}\right) \right] /n$ is precisely $B_{n}$ in
Theorem \ref{TH2}.
\end{remark}

\textbf{Proof : }

We have%
\begin{equation*}
\left\Vert T_{n}\right\Vert ^{2}=\frac{1}{n^{2}}\sum_{i=1}^{n}\left\Vert
\varepsilon _{i}\right\Vert ^{2}\left\langle \Gamma _{n}^{\dag
}X_{i},X_{n+1}\right\rangle ^{2}+\frac{1}{n^{2}}\sum_{i\neq i^{\prime
}}^{n}\left\langle \varepsilon _{i},\varepsilon _{i^{\prime }}\right\rangle
\left\langle \Gamma _{n}^{\dag }X_{i},X_{n+1}\right\rangle \left\langle
\Gamma _{n}^{\dag }X_{i^{\prime }},X_{n+1}.\right\rangle
\end{equation*}%
We take expectations in the display above and we note that the distribution
of each member of the first series on the right hand side does not depend on 
$n$ or $i$ and, due to linearity of expectation and $\mathbb{E}\left(
\varepsilon _{i}|X_{i}\right) =0$, the expectation of the second series is
null, hence%
\begin{align*}
\mathbb{E}\left\Vert T_{n}\right\Vert ^{2}& =\frac{1}{n}\mathbb{E}\left[
\left\Vert \varepsilon _{1}\right\Vert ^{2}\left\langle \Gamma _{n}^{\dag
}X_{1},X_{n+1}\right\rangle ^{2}\right] \\
& =\frac{1}{n}\mathbb{E}\left\{ \mathbb{E}\left[ \left\Vert \varepsilon
_{1}\right\Vert ^{2}\left\langle \Gamma _{n}^{\dag
}X_{1},X_{n+1}\right\rangle ^{2}|\varepsilon _{1},X_{1},...,X_{n}\right]
\right\} \\
& =\frac{1}{n}\mathbb{E}\left[ \left\Vert \varepsilon _{1}\right\Vert
^{2}|X_{1}\right] \mathbb{E}\left\langle \Gamma _{n}^{\dag }\Gamma \Gamma
_{n}^{\dag }X_{1},X_{1}\right\rangle .
\end{align*}%
We focus on $\mathbb{E}\left\langle \Gamma _{n}^{\dag }\Gamma \Gamma
_{n}^{\dag }X_{1},X_{1}\right\rangle $ and we see that this expectation is
nothing but the expectation of the trace of the operator $\Gamma _{n}^{\dag
}\Gamma \Gamma _{n}^{\dag }\cdot \left( X_{1}\otimes X_{1}\right) $, hence%
\begin{equation*}
\mathbb{E}\left\langle \Gamma _{n}^{\dag }\Gamma \Gamma _{n}^{\dag
}X_{1},X_{1}\right\rangle =\mathbb{E}\left\langle \Gamma _{n}^{\dag }\Gamma
\Gamma _{n}^{\dag }X_{i},X_{i}\right\rangle =\mathbb{E}\left[ \mathrm{tr}%
\Gamma _{n}^{\dag }\Gamma \Gamma _{n}^{\dag }\cdot \left( X_{i}\otimes
X_{i}\right) \right] ,
\end{equation*}%
and%
\begin{align*}
\mathbb{E}\left\langle \Gamma _{n}^{\dag }\Gamma \Gamma _{n}^{\dag
}X_{1},X_{1}\right\rangle & =\frac{1}{n}\mathbb{E}\left[ \mathrm{tr}\Gamma
_{n}^{\dag }\Gamma \Gamma _{n}^{\dag }\cdot \sum_{i=1}^{n}\left(
X_{i}\otimes X_{i}\right) \right] \\
& =\mathbb{E}\mathrm{tr}\left[ \Gamma _{n}^{\dag }\Gamma \Gamma _{n}^{\dag
}\Gamma _{n}\right] =\mathbb{E}\mathrm{tr}\left[ \Gamma _{n}^{\dag }\Gamma 
\widehat{\Pi }_{k}\right] =\mathbb{E}\mathrm{tr}\left[ \widehat{\Pi }%
_{k}\Gamma _{n}^{\dag }\Gamma \right] =\mathrm{tr}\left[ \Gamma \mathbb{E}%
\Gamma _{n}^{\dag }\right] .
\end{align*}%
At last we get :%
\begin{align*}
\mathbb{E}\left\langle \Gamma _{n}^{\dag }\Gamma \Gamma _{n}^{\dag
}X_{1},X_{1}\right\rangle & =\mathrm{tr}\left[ \Gamma \Gamma ^{\dag }\right]
+\mathrm{tr}\left[ \Gamma \mathbb{E}\left( \Gamma _{n}^{\dag }-\Gamma ^{\dag
}\right) \right] \\
& =k+\mathrm{tr}\left[ \Gamma \mathbb{E}\left( \Gamma _{n}^{\dag }-\Gamma
^{\dag }\right) \right] .
\end{align*}

>From Lemma \ref{approx-var} just below, we deduce that $\mathrm{tr}\left[
\Gamma \mathbb{E}\left( \Gamma _{n}^{\dag }-\Gamma ^{\dag }\right) \right]
=o\left( k\right) $, which finishes the proof of Proposition \ref{variance}.

\begin{lemma}
\label{approx-var}We have $\mathrm{tr}\left[ \Gamma \mathbb{E}\left( \Gamma
_{n}^{\dag }-\Gamma ^{\dag }\right) \right] \leq Ck^{2}\left( \log k\right)
/n,$ where $C$ does not depend on $S$, $n$ or $k.$ The preceding bound is an 
$o\left( k\right) $ since $k\left( \log k\right) /n\rightarrow 0$.
\end{lemma}

\textbf{Proof : }We focus on 
\begin{align*}
\left( \Gamma _{n}^{\dag }-\Gamma ^{\dag }\right) & =-\int_{C_{n}}\frac{1}{z}%
\left( zI-\Gamma _{n}\right) ^{-1}\left( \Gamma _{n}-\Gamma \right) \left(
zI-\Gamma \right) ^{-1}dz \\
& =-\int_{C_{n}}\frac{1}{z}\left( zI-\Gamma \right) ^{-1}\left( \Gamma
_{n}-\Gamma \right) \left( zI-\Gamma \right) ^{-1}dz \\
& -\int_{C_{n}}\frac{1}{z}\left( zI-\Gamma _{n}\right) ^{-1}\left( \Gamma
_{n}-\Gamma \right) \left( zI-\Gamma \right) ^{-1}\left( \Gamma _{n}-\Gamma
\right) \left( zI-\Gamma \right) ^{-1}dz.
\end{align*}

But $\mathbb{E}\int_{C_{n}}\frac{1}{z}\left( zI-\Gamma \right) ^{-1}\left(
\Gamma _{n}-\Gamma \right) \left( zI-\Gamma \right) ^{-1}dz=\int_{C_{n}}%
\frac{1}{z}\left( zI-\Gamma \right) ^{-1}\mathbb{E}\left( \Gamma _{n}-\Gamma
\right) \left( zI-\Gamma \right) ^{-1}dz=0$ so we consider the second term
above%
\begin{align*}
R_{n}& =\int_{C_{n}}\frac{1}{z}\left( zI-\Gamma _{n}\right) ^{-1}\left(
\Gamma _{n}-\Gamma \right) \left( zI-\Gamma \right) ^{-1}\left( \Gamma
_{n}-\Gamma \right) \left( zI-\Gamma \right) ^{-1}dz \\
& =\int_{C_{n}}\frac{1}{z}\left( zI-\Gamma \right) ^{-1/2}T_{n}\left(
z\right) A_{n}\left( z\right) A_{n}\left( z\right) \left( zI-\Gamma \right)
^{-1/2}dz,
\end{align*}%
where%
\begin{equation*}
T_{n}\left( z\right) =\left( zI-\Gamma \right) ^{1/2}\left( zI-\Gamma
_{n}\right) ^{-1}\left( zI-\Gamma \right) ^{1/2},\quad A_{n}\left( z\right)
=\left( zI-\Gamma \right) ^{-1/2}\left( \Gamma _{n}-\Gamma \right) \left(
zI-\Gamma \right) ^{-1/2}
\end{equation*}%
whence%
\begin{align*}
& \mathrm{tr}\left[ \Gamma R_{n}\right] =\sum_{j=1}^{+\infty }\int_{C_{n}}%
\frac{\lambda _{j}}{z-\lambda _{j}}\left\langle T_{n}\left( z\right)
A_{n}\left( z\right) A_{n}\left( z\right) \left( e_{j}\right) ,\left(
e_{j}\right) \right\rangle dz \\
& =\int_{C_{n}}\sum_{j=1}^{+\infty }\frac{\lambda _{j}}{z-\lambda _{j}}%
\left\langle T_{n}\left( z\right) A_{n}\left( z\right) A_{n}\left( z\right)
\left( e_{j}\right) ,\left( e_{j}\right) \right\rangle dz=\int_{C_{n}}%
\mathrm{tr}\left[ \left( zI-\Gamma \right) ^{-1}\Gamma T_{n}\left( z\right)
A_{n}\left( z\right) A_{n}\left( z\right) \right] dz,
\end{align*}%
and $\left\vert \mathrm{tr}\left[ \Gamma R_{n}\right] \right\vert \leq
\int_{C_{n}}\left[ \left\Vert \left( zI-\Gamma \right) ^{-1}\Gamma
T_{n}\left( z\right) \right\Vert _{\infty }\left\Vert A_{n}\left( z\right)
\right\Vert _{\mathcal{L}_{2}}^{2}\right] dz.$ Indeed, if we denote 
\begin{equation*}
\mathrm{tr}\left[ \left( zI-\Gamma \right) \Gamma T_{n}\left( z\right)
A_{n}\left( z\right) A_{n}\left( z\right) \right] =\mathrm{tr}\left[
A_{n}\left( z\right) \widetilde{T}_{n}\left( z\right) A_{n}\left( z\right) %
\right]
\end{equation*}%
with $\widetilde{T}_{n}\left( z\right) =\Gamma ^{1/2}\left( zI-\Gamma
_{n}\right) ^{-1}\Gamma ^{1/2}$ symmetric, we obtain 
\begin{equation*}
\mathrm{tr}\left[ A_{n}\left( z\right) \widetilde{T}_{n}\left( z\right)
A_{n}\left( z\right) \right] =\left\Vert \widetilde{T}_{n}^{1/2}\left(
z\right) A_{n}\left( z\right) \right\Vert _{\mathcal{L}_{2}}^{2}\leq
\left\Vert \widetilde{T}_{n}^{1/2}\left( z\right) \right\Vert _{\infty
}^{2}\left\Vert A_{n}\left( z\right) \right\Vert _{\mathcal{L}_{2}}^{2}.
\end{equation*}%
Now let us fix $m$. We have $\left\Vert \widetilde{T}_{n}^{1/2}\left(
z\right) \right\Vert _{\infty }^{2}\leq \left\Vert \widetilde{T}_{n}\left(
z\right) \right\Vert _{\infty }$ and $\sup_{z\in \mathcal{B}_{m}}\left\Vert 
\widetilde{T}_{n}\left( z\right) \right\Vert _{\infty }\leq Cm\quad a.s.$
The first inequality comes from the fact that $\widetilde{T}_{n}\left(
z\right) $ is symmetric, hence $\left\Vert \widetilde{T}_{n}\left( z\right)
\right\Vert _{\infty }=\sup_{\left\Vert u\right\Vert \leq 1}\left\vert
\left\langle \widetilde{T}_{n}\left( z\right) u,u\right\rangle \right\vert $%
. The last one comes from :%
\begin{equation*}
\widetilde{T}_{n}\left( z\right) =\Gamma ^{1/2}\left( zI-\Gamma \right)
^{-1/2}\left( zI-\Gamma \right) ^{1/2}\left( zI-\Gamma _{n}\right)
^{-1}\left( zI-\Gamma \right) ^{1/2}\left( zI-\Gamma \right) ^{-1/2}\Gamma
^{1/2},
\end{equation*}%
and 
\begin{equation*}
\left\Vert \widetilde{T}_{n}\left( z\right) \right\Vert _{\infty }\leq
\left\Vert \left( zI-\Gamma \right) ^{1/2}\left( zI-\Gamma _{n}\right)
^{-1}\left( zI-\Gamma \right) ^{1/2}\right\Vert _{\infty }\left\Vert \left(
zI-\Gamma \right) ^{-1}\Gamma \right\Vert _{\infty }.
\end{equation*}%
These facts prove (\ref{interm}). Now, by Lemma \ref{L1}, we can write $%
\mathbb{E}\left\Vert A_{n}\left( z\right) \right\Vert _{\mathcal{L}%
_{2}}^{2}\leq C\left( j\log j\right) ^{2}/n,$and consequently $\mathbb{E}%
\left\vert \mathrm{tr}\left[ \Gamma R_{n}\right] \right\vert \leq
C\sum_{j=1}^{k}\delta _{j}\frac{j^{3}\left( \log j\right) ^{2}}{n}=C\frac{1}{%
n}\sum_{j=1}^{k}\left( \lambda _{j}-\lambda _{j+1}\right) j^{3}\left( \log
j\right) ^{2}.$ By an Abel transform, we get :%
\begin{align*}
\sum_{j=1}^{k}\left( \lambda _{j}-\lambda _{j+1}\right) j^{3}\left( \log
j\right) ^{2}& \leq \frac{\lambda _{k+1}}{n}k^{3}\left( \log k\right) ^{2}+%
\frac{1}{n}\sum_{j=1}^{k}\lambda _{j}j^{2}\left( \log j\right) ^{2} \\
& \leq \frac{k^{2}\left( \log k\right) }{n}+\frac{1}{n}\sum_{j=1}^{k}j\left(
\log j\right) \leq \frac{k^{2}\left( \log k\right) }{n},
\end{align*}%
which yields $\mathbb{E}\left\vert \mathrm{tr}\left[ \Gamma R_{n}\right]
\right\vert \leq C\frac{k^{2}\left( \log k\right) }{n},$where $C$ is a
universal constant. Finally $\left\vert \mathrm{tr}\left[ \Gamma \mathbb{E}%
\left( \Gamma _{n}^{\dag }-\Gamma ^{\dag }\right) \right] \right\vert /k\leq
Ck\left( \log k\right) /n\rightarrow 0$ and we proved Lemma \ref{approx-var}%
. Now we are ready to turn to Theorem \ref{TH2}.\bigskip

\textbf{Proof of Theorem \ref{TH2} :}

>From equation (\ref{decomp-pred}), we obtain%
\begin{equation*}
\mathbb{E}\left\Vert S_{n}\left( X_{n+1}\right) -S\left( X_{n+1}\right)
\right\Vert ^{2}=\mathbb{E}\left\Vert S\widehat{\Pi }_{k}\left(
X_{n+1}\right) -S\left( X_{n+1}\right) \right\Vert ^{2}+\mathbb{E}\left\Vert 
\frac{1}{n}\sum_{i=1}^{n}\varepsilon _{i}\left\langle \Gamma _{n}^{\dag
}X_{i},X_{n+1}\right\rangle \right\Vert ^{2}.
\end{equation*}%
>From Proposition \ref{variance} followed by Lemma \ref{approx-var}, the
second term is $\frac{\sigma _{\varepsilon }^{2}}{n}k+B_{n}.$ It follows
from Proposition \ref{ks} and basic calculations that :%
\begin{equation*}
\mathbb{E}\left\Vert S\widehat{\Pi }_{k}\left( X_{n+1}\right) -S\left(
X_{n+1}\right) \right\Vert ^{2}=\mathbb{E}\left\Vert S\left( \Pi
_{k}-I\right) \left( X_{n+1}\right) \right\Vert ^{2}+A_{n},
\end{equation*}%
where $A_{n}$ matches the bound of the Theorem. At last $\mathbb{E}%
\left\Vert S\left( \Pi _{k}-I\right) \left( X_{n+1}\right) \right\Vert
^{2}=\sum_{j\geq k+1}\lambda _{j}\left\Vert Se_{j}\right\Vert ^{2}$ which
finishes the proof.\bigskip

\textbf{Proof of Theorem \ref{TH2bis} :}

Our proof follows the lines of Cardot, Johannes (2010) through a modified
version of Assouad's lemma.

To simplify notations we set $k_{n}^{\ast }=k_{n}.$ Take $S^{\theta
}=\sum_{j=1}^{k_{n}}\eta _{i}\omega _{i}e_{i}\otimes e_{1}$ where $\omega
_{i}\in \left\{ -1,1\right\} $ and $\theta =\left[ \omega _{1},...,\omega
_{k}\right] $ and $\eta _{i}\in \mathbb{R}^{+}$ will be fixed later such
that $S^{\theta }\in \mathcal{L}_{2}\left( \varphi ,C\right) $ for all $%
\theta $. Denote $\theta _{-i}=\left[ \omega _{1},...,-\omega
_{i},...,\omega _{k}\right] $ and $\mathbb{P}_{\theta }:=\mathbb{P}_{\theta }%
\left[ \left( Y_{1},X_{1}\right) ,...,\left( Y_{n},X_{n}\right) \right] $
denote the distribution of the data when $S=S^{\theta }$. Let $\rho $ stand
for Hellinger's affinity, $\rho \left( \mathbb{P}_{0},\mathbb{P}_{1}\right)
=\int \sqrt{d\mathbb{P}_{0}d\mathbb{P}_{1}}$ and $\mathbf{KL}\left( \mathbb{P%
}_{0},\mathbb{P}_{1}\right) $ for K\"{u}llback-Leibler divergence then $\rho
\left( \mathbb{P}_{0},\mathbb{P}_{1}\right) \geq \left( 1-\frac{1}{2}\mathbf{%
KL}\left( \mathbb{P}_{0},\mathbb{P}_{1}\right) \right) .$

Note that considering models based on $S^{\theta }$ above comes down to
projecting the model on a one-dimensional space. We are then faced with a
linear model with real output and finally confine ourselves to proving that
the optimal rate is unchanged (see Hall, Horowitz (2007)).

\begin{eqnarray*}
\mathcal{R}_{n}\left( T_{n}\right) &=&\sup_{S\in \mathcal{L}_{2}\left(
\varphi ,C\right) }\mathbb{E}\left\Vert \left( T_{n}-S\right) \Gamma
^{1/2}\right\Vert _{2}^{2}\geq \frac{1}{2^{k}}\sum_{\omega \in \left\{
-1,1\right\} ^{k}}\sum_{i=1}^{k_{n}}\lambda _{i}\mathbb{E}_{\theta
}\left\langle \left( T_{n}-S^{\theta }\right) e_{i},e_{1}\right\rangle ^{2}
\\
&=&\frac{1}{2^{k}}\sum_{\omega \in \left\{ -1,1\right\} ^{k}}\frac{1}{2}%
\sum_{i=1}^{k_{n}}\lambda _{i}\left[ \mathbb{E}_{\theta }\left\langle \left(
T_{n}-S^{\theta }\right) e_{i},e_{1}\right\rangle ^{2}+\mathbb{E}_{\theta
_{_{-i}}}\left\langle \left( T_{n}-S^{_{\theta _{_{-i}}}}\right)
e_{i},e_{1}\right\rangle ^{2}\right] \\
&\geq &\frac{1}{2^{k}}\sum_{\omega \in \left\{ -1,1\right\}
^{k}}\sum_{i=1}^{k_{n}}\lambda _{i}\eta _{i}^{2}\rho ^{2}\left( \mathbb{P}%
_{\theta },\mathbb{P}_{\theta _{-i}}\right)
\end{eqnarray*}

The last line was obtained by a slight variant of the bound (A.9) in Cardot,
Johannes (2010), p.405 detailed below :%
\begin{eqnarray*}
\rho \left( \mathbb{P}_{\theta },\mathbb{P}_{\theta _{-i}}\right) &\leq
&\int \frac{\left\langle \left( T_{n}-S^{\theta }\right)
e_{i},e_{1}\right\rangle }{\left\vert \left\langle \left( S^{\theta
_{-i}}-S^{\theta }\right) e_{i},e_{1}\right\rangle \right\vert }\sqrt{d%
\mathbb{P}_{0}d\mathbb{P}_{1}}+\int \frac{\left\langle \left(
T_{n}-S^{\theta _{-i}}\right) e_{i},e_{1}\right\rangle }{\left\vert
\left\langle \left( S^{\theta _{-i}}-S^{\theta }\right)
e_{i},e_{1}\right\rangle \right\vert }\sqrt{d\mathbb{P}_{0}d\mathbb{P}_{1}}
\\
&\leq &\frac{1}{2\eta _{i}}\left( \int \left\langle \left( T_{n}-S^{\theta
}\right) e_{i},e_{1}\right\rangle ^{2}d\mathbb{P}_{\theta }\right)
^{1/2}+\left( \int \left\langle \left( T_{n}-S^{\theta _{-i}}\right)
e_{i},e_{1}\right\rangle \mathbb{P}_{\theta _{-i}}\right) ^{1/2}
\end{eqnarray*}

by Cauchy-Schwartz inequality and since $\left\vert \left\langle \left(
S^{_{\theta _{-i}}}-S^{\theta }\right) e_{i},e_{1}\right\rangle \right\vert
=2\eta _{i}$. Then 
\begin{equation*}
2\eta _{i}^{2}\rho ^{2}\left( \mathbb{P}_{\theta },\mathbb{P}_{\theta
_{-i}}\right) \leq \mathbb{E}_{\theta }\left\langle \left( T_{n}-S^{\theta
}\right) e_{i},e_{1}\right\rangle ^{2}+\mathbb{E}_{\theta _{-i}}\left\langle
\left( T_{n}-S^{\theta _{-i}}\right) e_{i},e_{1}\right\rangle ^{2}
\end{equation*}%
yields :%
\begin{equation*}
\mathcal{R}_{n}\left( T_{n}\right) \geq \inf_{\omega \in \left\{
-1,1\right\} ^{k}}\inf_{i}\rho \left( \mathbb{P}_{\theta },\mathbb{P}%
_{\theta _{-i}}\right) \sum_{i}\lambda _{i}\eta _{i}^{2}
\end{equation*}

We show below that $\mathbf{KL}\left( \mathbb{P}_{\theta },\mathbb{P}%
_{\theta _{-i}}\right) \leq 4n\lambda _{i}\eta _{i}^{2}/\sigma _{1}^{2}$.
Choosing $\eta _{i}=\sigma _{1}/2\sqrt{n\lambda _{i}}$ for $1\leq i\leq
k_{n} $ gives $S^{\theta }\in \mathcal{L}_{2}\left( \varphi ,1\right) $ and $%
\sup_{\omega ,i}\mathbf{KL}\left( \mathbb{P}_{\theta },\mathbb{P}_{\theta
_{-i}}\right) \leq 1$, $\inf_{\omega ,i}\rho \left( \mathbb{P}_{\theta },%
\mathbb{P}_{\theta _{-i}}\right) \leq 1/2$ and%
\begin{equation*}
\mathcal{R}_{n}\left( T_{n}\right) \geq \frac{1}{2}\sum_{i=1}^{k_{n}}\lambda
_{i}\eta _{i}^{2}=\frac{1}{2}\frac{k_{n}}{n}
\end{equation*}%
whatever the choise of the estimate $T_{n}$. This proves the lower bound :%
\begin{equation*}
\lim \sup_{n\rightarrow +\infty }\varphi _{n}^{-1}\inf_{T_{n}}\sup_{S\in 
\mathcal{L}_{2}\left( \varphi ,C\right) }E\left\Vert \left( T_{n}-S\right)
\Gamma ^{1/2}\right\Vert ^{2}>\frac{1}{2},
\end{equation*}%
and the Theorem stems from this last display.

We finish by proving that $\mathbf{KL}\left( \mathbb{P}_{\theta },\mathbb{P}%
_{\theta _{-i}}\right) \leq 4n\lambda _{i}\eta _{i}^{2}/\sigma _{1}^{2}.$ It
suffices to notice that 
\begin{equation*}
\mathbf{KL}\left( \mathbb{P}_{\theta },\mathbb{P}_{\theta _{-i}}\right)
=\int \log \left( d\mathbb{P}_{\theta |X}/d\mathbb{P}_{\theta
_{-i}|X}\right) d\mathbb{P}_{\theta }
\end{equation*}%
where $\mathbb{P}_{\theta |X}$ stand for the likelihood of $Y$ conidtionally
to $X$. In this Hilbert setting we must clarify the existence of this
likelihood ratio. It suffices to prove that $\mathbb{P}_{\theta |X}\left(
Y\right) \ll \mathbb{P}_{0|X}\left( Y\right) $ which in turn is true when $%
S^{\theta }X$ belongs to the RKHS associated to $\varepsilon $ (see Lifshits
(1995)). With other words we need that almost surely $\Gamma _{\varepsilon
}^{-1/2}S^{\theta }X$ is finite where $\Gamma _{\varepsilon }$ is the
covariance operator of the noise. But $\Gamma _{\varepsilon
}^{-1/2}S^{\theta }=S^{\theta }/\sigma _{1}$. Set $\omega _{l}^{\prime
}=\omega _{l}$ if $l\neq i$ with $\omega _{i}^{\prime }=-\omega _{i}$ :%
\begin{eqnarray*}
\log \frac{d\mathbb{P}_{\theta |X}\left( Y\right) }{d\mathbb{P}_{\theta
_{-i}|X}\left( Y\right) } &=&-\left( \left\langle Y,e_{1}\right\rangle
-\sum_{l=1}^{k_{n}}\omega _{l}\eta _{l}\left\langle X,e_{l}\right\rangle
\right) ^{2}+\left( \left\langle Y,e_{1}\right\rangle
-\sum_{l=1}^{k_{n}}\omega _{l}^{\prime }\eta _{l}\left\langle
X,e_{l}\right\rangle \right) ^{2} \\
&=&-2\omega _{i}\eta _{i}\frac{\left\langle X,e_{i}\right\rangle }{\sigma
_{1}^{2}}\left( 2\left\langle \varepsilon ,e_{1}\right\rangle
+\sum_{l=1}^{k_{n}}\omega _{l}\eta _{l}\left\langle X,e_{l}\right\rangle
-\sum_{l=1}^{k_{n}}\omega _{l}^{\prime }\eta _{l}\left\langle
X,e_{l}\right\rangle \right)  \\
&=&-2\omega _{i}\eta _{i}\frac{\left\langle X,e_{i}\right\rangle }{\sigma
_{1}^{2}}\left( 2\left\langle \varepsilon ,e_{1}\right\rangle +2\omega
_{i}\eta _{i}\left\langle X,e_{i}\right\rangle \right) 
\end{eqnarray*}%
and $\mathbb{E}_{\theta }\left[ \log d\mathbb{P}_{\theta |X}\left( Y\right)
/d\mathbb{P}_{\theta _{-i}|X}\left( Y\right) \right] =4\eta _{i}^{2}\mathbb{E%
}_{\theta }\left\langle X,e_{i}\right\rangle ^{2}/\sigma _{1}^{2}=4\eta
_{i}^{2}\mathbb{\lambda }_{i}/\sigma _{1}^{2}$\bigskip 

Now we focus on the problem of weak convergence.\bigskip

\textbf{Proof of Theorem \ref{TH1} :}

Consider (\ref{decomp}). We claim that weak convergence of $S_{n}$ will
depend on the series $\left( 1/n\right) \sum_{i=1}^{n}\varepsilon
_{i}\otimes \Gamma _{n}^{\dag }X_{i}$. This fact can be checked by
inspecting the proof of Theorem \ref{TH2}. We are going to prove that $%
\left( 1/n\right) \sum_{i=1}^{n}\varepsilon _{i}\otimes \Gamma ^{\dag }X_{i}$
cannot converge for the classical (supremum) operator norm. We replace the
random $\Gamma _{n}^{\dag }$ by the non-random $\Gamma ^{\dag }$. It is
plain that non-convergence of the second series implies non-convergence of
the first. Suppose that for some sequence $\alpha _{n}\uparrow +\infty $ the
centered series $\left( \alpha _{n}/n\right) \sum_{i=1}^{n}\varepsilon
_{i}\otimes \Gamma ^{\dag }X_{i}\overset{w}{\rightarrow }Z,$in operator
norm, where $Z$ is a fixed random operator (not necessarily gaussian). Then
for all fixed $x$ and $y$ in $H,\frac{\alpha _{n}}{n}\sum_{i=1}^{n}\left%
\langle \varepsilon _{i},y\right\rangle \left\langle \Gamma ^{\dag
}X_{i},x\right\rangle \overset{w}{\rightarrow }\left\langle
Zx,y\right\rangle ,$as real random variables. First take $x$ in the domain
of $\Gamma ^{-1}$. From $\left\Vert \Gamma ^{-1}x\right\Vert <+\infty $, we
see that $\mathbb{E}\left\langle \varepsilon _{i},y\right\rangle
^{2}\left\langle \Gamma ^{\dag }X_{i},x\right\rangle ^{2}<+\infty $ implies
that $\alpha _{n}=\sqrt{n}$ (and $Z$ is gaussian since we apply the central
limit theorem for independent random variables). Now take a $x$ such that $%
\left\Vert \Gamma ^{-1}x\right\Vert =+\infty $, then $\mathbb{E}\left\langle
\varepsilon _{1},y\right\rangle ^{2}\left\langle \Gamma ^{\dag
}X_{1},x\right\rangle ^{2}=\mathbb{E}\left\langle \varepsilon
_{1},y\right\rangle ^{2}\mathbb{E}\left\langle \Gamma ^{\dag
}x,x\right\rangle $, and is is easily seen from the definition of $\Gamma
^{\dag }$ that $\mathbb{E}\left\langle \Gamma ^{\dag }x,x\right\rangle $
-which is positive and implicitely depend on $n$ through $k$- tends to
infinity. Consequently $\left( 1/\sqrt{n}\right) \sum_{i=1}^{n}\varepsilon
_{i}\otimes \Gamma ^{\dag }X_{i}$ cannot converge weakly anymore since the
margins related to the $x$'s do not converge in distribution. This proves
the Theorem.\bigskip

The two next Lemmas prepare the proof of Theorem \ref{TH3}. We set $T_{n}=%
\frac{1}{n}\sum_{i=1}^{n}\varepsilon _{i}\left\langle \Gamma _{n}^{\dag
}X_{i},X_{n+1}\right\rangle $ and this series is the crucial term that
determines weak convergence. We go quickly through the first Lemma since it
is close to Lemma 8 p.355 in Cardot, Mas, Sarda (2007).

\begin{lemma}
\label{conv.loi.proj.pred}Fix $x$ in $H,$ then $\sqrt{n/k_{n}}\left\langle
T_{n},x\right\rangle \overset{w}{\rightarrow }\mathcal{N}\left( 0,\sigma
_{\varepsilon ,x}^{2}\right) $, where $\sigma _{\varepsilon ,x}^{2}=\mathbb{E%
}\left\langle \varepsilon _{k},x\right\rangle ^{2}$.
\end{lemma}

\textbf{Proof :} Let $\mathcal{F}_{n}$ be the $\sigma $-algebra generated by 
$\left( \varepsilon _{1},...,\varepsilon _{n},X_{1},...,X_{n}\right) $. We
see that $Z_{i,n}^{x}=\left\langle \varepsilon _{i},x\right\rangle
\left\langle \Gamma _{n}^{\dag }X_{i},X_{n+1}\right\rangle $ is a
real-valued martingale difference, besides%
\begin{equation*}
\mathbb{E}\left[ \left( Z_{i,n}^{x}\right) ^{2}|\mathcal{F}_{n}\right]
=\sigma _{\varepsilon ,x}^{2}\left\langle \Gamma _{n}^{\dag
}X_{i},X_{n+1}\right\rangle ^{2}.
\end{equation*}%
Applying Lemma \ref{approx-var} and results by McLeish (1974) on weak
convergence for martingale differences arrays yields the Lemma.

\begin{lemma}
\label{flat-conc}The random sequence $\sqrt{\frac{k_{n}}{n}}T_{n}$ is flatly
concentrated and uniformly tight. In fact, if $\mathcal{P}_{m}$ is the
projection operator on the $m$ first eigenvectors of $\Gamma _{\varepsilon }$
and $\eta >0$ is a real number 
\begin{equation*}
\limsup_{m\rightarrow +\infty }\sup_{n}\mathbb{P}\left( \left\Vert \sqrt{%
\frac{n}{k_{n}}}\left( I-\mathcal{P}_{m}\right) T_{n}\right\Vert >\eta
\right) =0.
\end{equation*}
\end{lemma}

\textbf{Proof :} Let $\mathcal{P}_{m}$ be the projection operator on the $m$
first eigenvectors of $\Gamma _{\varepsilon }$. For $\sqrt{k_{n}/n}T_{n}$ to
be flatly concentrated it is sufficient to prove that for any $\eta >0$, 
\begin{equation*}
\limsup_{m\rightarrow +\infty }\sup_{n}\mathbb{P}\left( \left\Vert \sqrt{%
\frac{n}{k_{n}}}\left( I-\mathcal{P}_{m}\right) T_{n}\right\Vert >\eta
\right) =0.
\end{equation*}%
We have :%
\begin{align*}
& \mathbb{P}\left( \left\Vert \sqrt{\frac{n}{k_{n}}}\left( I-\mathcal{P}%
_{m}\right) T_{n}\right\Vert >\eta \right) \\
& \leq \frac{1}{\eta ^{2}}\mathbb{E}\left\Vert \sqrt{\frac{n}{k_{n}}}\left(
I-\mathcal{P}_{m}\right) T_{n}\right\Vert ^{2}=\frac{1}{\eta ^{2}k_{n}}%
\mathbb{E}\left\langle \Gamma _{n}^{\dag }X_{1},X_{n+1}\right\rangle ^{2}%
\mathbb{E}\left\Vert \left( I-\mathcal{P}_{m}\right) \varepsilon
_{1}\right\Vert ^{2}.
\end{align*}%
We see first that $\sup_{n}\mathbb{P}\left( \left\Vert \sqrt{\frac{n}{k_{n}}}%
\left( I-\mathcal{P}_{m}\right) T_{n}\right\Vert >\eta \right) \leq \frac{C}{%
\eta ^{2}}\mathbb{E}\left\Vert \left( I-\mathcal{P}_{m}\right) \varepsilon
_{1}\right\Vert ^{2}$ where $C$ is some constant and once again following
Lemma \ref{approx-var}. Now it is plain that 
\begin{equation*}
\limsup_{m\rightarrow +\infty }\mathbb{E}\left\Vert \left( I-\mathcal{P}%
_{m}\right) \varepsilon _{1}\right\Vert ^{2}=0,
\end{equation*}%
because $\mathcal{P}_{m}$ was precisely chosen to be projector on the $m$
first eigenvectors of the trace-class operator $\Gamma _{\varepsilon }$. In
fact $\mathbb{E}\left\Vert \left( I-\mathcal{P}_{m}\right) \varepsilon
_{1}\right\Vert ^{2}=\mathrm{tr}\left[ \left( I-\mathcal{P}_{m}\right)
\Gamma _{\varepsilon }\left( I-\mathcal{P}_{m}\right) \right] ,$and this
trace is nothing but the series summing the eigenvalues of $\Gamma
_{\varepsilon }$ from order $m+1$ to infinity, hence the result.\bigskip

\textbf{Proof of Theorem \ref{TH3} : }We only prove the second part of the
theorem : weak convergence with no bias. The first part follows immediately.
We start again from the decomposition (\ref{decomp-pred}). As announced just
above, the two first terms vanish with respect to convergence in
distribution. For $S\left[ \widehat{\Pi }_{k}-\Pi _{k}\right] \left(
X_{n+1}\right) $, we invoke Proposition \ref{ks} to claim that, whenever $%
k^{2}\log ^{2}k/n\rightarrow 0$, then $\left( n/k\right) \mathbb{E}%
\left\Vert S\left[ \widehat{\Pi }_{k}-\Pi _{k}\right] \left( X_{n+1}\right)
\right\Vert ^{2}\rightarrow 0$ and we just have to deal with the first term,
related to bias : $S\left( \Pi _{k}-I\right) \left( X_{n+1}\right) .$

Assume first that the mean square of the latter reminder, $\left( n/k\right)
\sum_{j=k+1}^{+\infty }\lambda _{j}\left\Vert S\left( e_{j}\right)
\right\Vert ^{2}$, decays to zero. Then the proof of the Theorem is
immediate from Lemmas \ref{conv.loi.proj.pred} and \ref{flat-conc}. The
sequence $\sqrt{n/k_{n}}T_{n}$ is uniformly tight and its finite dimensional
distributions (in the sense of "all finite-dimensional projections of $\sqrt{%
n/k_{n}}T_{n}$") converge weakly to $\mathcal{N}\left( 0,\sigma
_{\varepsilon ,x}^{2}\right) $. This is enough to claim that Theorem \ref%
{TH3} holds. We refer for instance to de Acosta (1970) or Araujo and Gin\'{e}
(1980) for checking the validity of this conclusion.

Finally, the only fact to be proved is $\lim_{n\rightarrow +\infty }\left(
n/k\right) \sum_{j=k+1}^{+\infty }\lambda _{j}\left\Vert S\left(
e_{j}\right) \right\Vert ^{2}=0$ when tightening conditions on the sequence $%
k_{n}.$ This looks like an Abelian theorem which could be proved by special
techniques but we prove it in a simple direct way. First, we know by
previous remarks (since $\lambda _{j}$ and $\left\Vert S\left( e_{j}\right)
\right\Vert ^{2}$ are convergent series) that $\lambda _{j}\left\Vert
S\left( e_{j}\right) \right\Vert ^{2}=\tau _{j}\left( j^{2}\log ^{2}j\right)
,$where $\tau _{j}$ tends to zero. Taking as in the first part of the
theorem $n=k^{2}\log ^{2}k/\sqrt{\gamma _{k}}$, we can focus on $%
\lim_{k+\infty }\frac{k\log ^{2}k}{\sqrt{\gamma _{k}}}\sum_{j=k+1}^{+\infty
}\tau _{j}/\left( j^{2}\log ^{2}j\right) .$ We know that for a sufficiently
large $k$ and for all $j\geq k,$ $0\leq \tau _{j}\leq \epsilon $ where $%
\epsilon >0$ is fixed. Then%
\begin{align*}
\frac{1}{\sqrt{\gamma _{k}}}\sum_{j=k+1}^{+\infty }\tau _{j}\frac{k\log ^{2}k%
}{j^{2}\log ^{2}j}& =\frac{1}{\sqrt{\gamma _{k}}}\sum_{m=1}^{+\infty
}\sum_{j=km+1}^{km+k}\tau _{j}\frac{k\log ^{2}k}{j^{2}\log ^{2}j} \\
& \leq \frac{1}{\sqrt{\gamma _{k}}}\sum_{m=1}^{+\infty }\left(
\sup_{km+1\leq j\leq km}\tau _{j}\right) \frac{k^{2}\log ^{2}k}{%
k^{2}m^{2}\log ^{2}km} \\
& \leq \frac{1}{\sqrt{\gamma _{k}}}\left( \sup_{k\leq j}\tau _{j}\right)
\sum_{m=1}^{+\infty }\frac{1}{m^{2}}=C\sqrt{\gamma _{k}}\rightarrow 0,
\end{align*}%
which removes the bias term and is the desired result.

\end{document}